\documentclass[review,hidelinks,onefignum,onetabnum]{siamart220329}
\usepackage{amsmath}
\usepackage{lipsum}
\usepackage{amsfonts}
\usepackage{graphicx}
\usepackage{epstopdf}
\usepackage{algorithmic}
\usepackage{amssymb}
\usepackage{amsopn}

\nolinenumbers

\ifpdf
\DeclareGraphicsExtensions{.eps,.pdf,.png,.jpg}
\else
\DeclareGraphicsExtensions{.eps}
\fi

\newsiamremark{remark}{Remark}
\newsiamremark{hypothesis}{Hypothesis}
\crefname{hypothesis}{Hypothesis}{Hypotheses}
\newsiamthm{claim}{Claim}
\headers{Convexification for a Coefficient Inverse Problem}
{M. V. Klibanov, J. Li, and Z. Yang}

\ifpdf
\hypersetup{
  pdftitle={Convexification for a Coefficient Inverse Problem},
  pdfauthor={M. V. Klibanov, J. Li, and Z. Yang}
}
\fi

\begin{document}

\title{Convexification for a Coefficient Inverse Problem for a System of Two
Coupled Nonlinear Parabolic Equations }
\author{Michael V. Klibanov \thanks{%
Department of Mathematics and Statistics, University of North Carolina at
Charlotte, Charlotte, NC, 28223, USA (mklibanv@charlotte.edu)} \and Jingzhi
Li \thanks{%
Department of Mathematics \& National Center for Applied Mathematics
Shenzhen \& SUSTech International Center for Mathematics, Southern
University of Science and Technology, Shenzhen 518055, P.~R.~China
(li.jz@sustech.edu.cn)} \and Zhipeng Yang \thanks{%
School of Mathematics and Statistics, Lanzhou University, Lanzhou 730000, P.
R. China (yangzp@lzu.edu.cn)} }
\maketitle

\begin{abstract}
A system of two coupled nonlinear parabolic partial differential equations
with two opposite directions of time is considered. In fact, this is the
so-called \textquotedblleft Mean Field Games System" (MFGS), which is
derived in the mean field games (MFG) theory. This theory has numerous
applications in social sciences. The topic of Coefficient Inverse Problems
(CIPs) in the MFG theory is in its infant age, both in theory and
computations. A numerical method for this CIP is developed. Convergence
analysis ensures the global convergence of this method. Numerical
experiments are presented.
\end{abstract}



\section{Introduction}

\label{sec:1}

A system of two nonlinear parabolic partial differential equations with two
opposite directions of time is considered. This is the so-called
\textquotedblleft Mean Field Games System" (MFGS), which is derived in the
mean field games (MFG) theory \cite{A}. This theory is about the mass
behavior of infinitely many agents, who are assumed to be rational. Due to a
rapidly growing number of applications of the MFG theory in social sciences,
it makes sense to study a variety of mathematical questions for the MFGS. We
are interested in a Coefficient Inverse Problem (CIP) for this system. Our
unknown coefficient $k\left( x\right) $ characterizes a reaction of a
controlled agent to an action applied at the point $x\in \Omega ,$ where $%
\Omega \subset \mathbb{R}^{n}$ is a bounded domain of interest, where the
agents are located \cite[page 632]{MFG1}. A significant complicating factor
of our CIP is that the unknown coefficient is involved in the MFGS along
with its first derivatives. Nevertheless, we still use the input data
resulting from only a single measurement event. H\"{o}lder stability
estimates and uniqueness for our CIP were obtained in \cite{MFG6}. In this
paper, we develop a version of the so-called convexification numerical
method for this problem. First, we construct this method. Next, we carry out
convergence analysis. This analysis ensures the global convergence property
of our technique. We end up with numerical experiments.

In simple terms, the notion of the global convergence means that the
requirement of the availability of a good first guess about the solution is
not imposed. To be more precise, we introduce the following Definition:

\textbf{Definition.} \emph{We call a numerical method for a CIP globally
convergent if there is a theorem, which guarantees that, as long as the
noise in the data tends to zero, iterates of this method converge to the
true solution of that CIP\ without any advanced knowledge of a sufficiently
small neighborhood of this solution. }

Since its initial introduction in 2006 in seminal publications of Lasry and
Lions \cite{LL1,LL2} and Huang, Caines, and Malham\'{e} \cite{Huang2,Huang1}%
, the MFG theory has found a number of applications in descriptions of
complex social phenomena. Among those applications, we mention, e.g. finance 
\cite{A,Trusov}, sociology \cite{Bauso}, election dynamics \cite{Chow},
interactions of electrical vehicles \cite{Co}, combating corruption \cite%
{Kol} and deep learning \cite{W}.

Since this theory is young, then the topic of CIPs for the MFGS is currently
in its infant stage. We now list publications on this topic. In the case of
data resulting from a single measurement event, such as ones considered
here, stability and uniqueness results were obtained in \cite{ImLY,MFG4,MFG6}%
, using the method of \cite{BukhKlib}. In \cite{Liu,Ren} uniqueness results
were obtained for the case of the data resulting from multiple measurements.
As to the numerical studies of CIPs for the MFGS, we refer to \cite%
{Chow,Ding}. The methods of \cite{Chow,Ding} are significantly different
from the one of the current paper. A numerical method, which is close to the
one of this paper, was developed by this research team in \cite{MFG8}.\ In 
\cite{MFG8} a version of the convexification method is developed for the
case when the coefficient $b\left( x\right) $ is unknown in the so-called
global interaction term%
\begin{equation}
\int\limits_{\Omega }K\left( x,y\right) m\left( y,t\right) dy=b\left(
x\right) \int\limits_{\Omega }\overline{K}\left( x,y\right) m\left(
y,t\right) dy,  \label{1.1}
\end{equation}%
where the part $\overline{K}\left( x,y\right) $ of the kernel $K\left(
x,y\right) $ is known, see section 2 for explanations of notations. However,
unlike the current paper, derivatives of the unknown coefficient $b\left(
x\right) $ are not involved in the MFGS. This means that the CIP of \cite%
{MFG8} is significantly simpler than the one considered here.

The development of the convexification method is caused by the desire to
avoid the appearance of local minima of conventional least squares cost
functionals for CIPs, see, e.g. \cite{B,BL,Giorgi,Gonch,MB} for such
functionals. Local minima cause the local convergence of gradient-like
minimization methods for those functionals. In other words, the convergence
of such a method to the true solution of a CIP can be rigorously guaranteed
only if the starting point of iterations is located in a sufficiently small
neighborhood of that solution. Unlike this, the convexification is a
globally convergent method in terms of the above Definition. The theory of
the convexification method was first published in \cite{Klib97}. More
recently a number of works about the convexification, combining the theory
with numerical studies, has been published, see, e.g. \cite%
{Bak,Baud1,Baud2,Baud3,Kpar,KL,MFG8} and references cited therein.

The convexification can be considered as a numerical development of the idea
of the paper \cite{BukhKlib}, which was originally concerned only with the
question of uniqueness theorems for multidimensional CIPs with single
measurement data. The method of \cite{BukhKlib} is based on Carleman
estimates, see, e.g. \cite{GY,Isakov,Ksurvey,KL,Lay,Ma,Rakesh} and
references cited therein for some follow up publications. These estimates
also play an important role in the convexification method. Carleman
estimates were first introduced in the MFG theory in \cite{MFG1} at the time
when the first version of this paper was posted as a prerpint on
www.arxiv.org, arXiv: 2302.10709v1, February 21, 2023, and they were used
since then in \cite{ImLY,MFG4,MFG6,MFG8}.

All functions considered below are real valued ones. In section 2 we bring
in the MFGS and state our CIP then. In particular, we indicate in this
section the main challenges in working with the MFGS. In section 3 we
describe a certain transformation procedure, which is the first step towards
the convexification. In section 4 we construct a weighted cost functional
for the transformed problem. This functional is the key to the
convexification. The weight is the so-called Carleman Weight Function (CWF),
i.e. the function, which is used as the weight in Carleman estimates for the
parabolic Partial Differential Operators involved in the transformed
problem. Then the focus is on the minimization of that functional. In
section 5 we formulate theorems of our convergence analysis. In section 6 we
prove the central theorem formulated in section 5. In section 7 we confirm
our theory by numerical experiments.

\section{Statement of the Coefficient Inverse Problem}

\label{sec:2}

\subsection{Statement of the CIP}

\label{sec:2.1}

Below $x=\left( x_{1},x_{2},...,x_{n}\right) \in \mathbb{R}^{n}$, $\overline{%
x}=(x_{2}$,$...,x_{n})$. To simplify the presentation, we assume that our
domain of interest $\Omega \subset \mathbb{R}^{n}$ is a rectangular prism.
Consider some numbers $a,b,A_{i},\gamma ,T>0,i=1,...,n$, where $a<b$ and $%
\gamma \in \left( 0,1\right) .$ Denote%
\begin{equation}
\left. 
\begin{array}{c}
\Omega =\left\{ x:a<x_{1}<b,-A_{i}<x_{i}<A_{i},i=2,...,n\right\} ,\text{ }
\\ 
\Omega _{1}=\left\{ y:-A_{i}<x_{i}<A_{i},i=2,...,n\right\} , \\ 
Q_{T}=\Omega \times \left( 0,T\right) ,\text{ }S_{T}=\partial \Omega \times
\left( 0,T\right) , \\ 
Q_{\gamma ,T}=\Omega \times \left( \left( 1-\gamma \right) T/2,\left(
1+\gamma \right) T/2\right) , \\ 
\Psi _{1}^{+}=\left\{ x\in \partial \Omega :x_{1}=b\right\} ,\text{ }\Psi
_{1}^{-}=\left\{ x\in \partial \Omega :x_{1}=a\right\} , \\ 
\Psi _{1,T}^{+}=\Psi _{1}^{+}\times \left( 0,T\right) ,\text{ }\Psi
_{1,T}^{-}=\Psi _{1}^{-}\times \left( 0,T\right) ,\Psi _{1,T}^{\pm }=\Psi
_{1,T}^{+}\cup \Psi _{1,T}^{-}, \\ 
\Psi _{i}^{\pm }=\left\{ x\in \partial \Omega :x_{i}=\pm A_{i}\right\} ,%
\text{ }\Psi _{iT}^{\pm }=\Psi _{i}^{\pm }\times \left( 0,T\right) ,\text{ }%
i=2,...,n.%
\end{array}%
\right.  \label{2.1}
\end{equation}%
Based on (\ref{2.1}),\emph{\ }denote 
\begin{equation*}
\left. 
\begin{array}{c}
H^{2,1}\left( Q_{\gamma ,T}\right) =\left\{ u:\left\Vert u\right\Vert
_{H^{2,1}\left( Q_{\gamma ,T}\right) }^{2}=\sum\limits_{\left\vert \alpha
\right\vert +2m\leq 2}\left\Vert D_{x}^{\alpha }\partial
_{t}^{m}u\right\Vert _{L_{2}\left( Q_{\gamma ,T}\right) }^{2}<\infty
\right\} , \\ 
H^{2,1}\left( \Psi _{iT}^{\pm }\right) =\left\{ 
\begin{array}{c}
u:\left\Vert u\right\Vert _{H^{2,1}\left( \Psi _{iT}^{\pm }\right)
}^{2}=\sum\limits_{j=1,j\neq i}^{n}\left\Vert u_{x_{j}}\right\Vert
_{L_{2}\left( \Psi _{iT}^{\pm }\right) }^{2}+ \\ 
+\sum\limits_{j,s=1,\left( j,s\right) \neq \left( i,i\right) }^{n}\left\Vert
u_{x_{j}x_{s}}\right\Vert _{L_{2}\left( \Psi _{iT}^{\pm }\right)
}^{2}+\sum\limits_{j=0}^{1}\left\Vert \partial _{t}^{j}u\right\Vert
_{L_{2}\left( \Psi _{iT}^{\pm }\right) }^{2}<\infty%
\end{array}%
\right\} , \\ 
H^{2,1}\left( S_{T}\right) =\left\{ u:\left\Vert u\right\Vert
_{H^{2,1}\left( S_{T}\right) }^{2}=\sum\limits_{i=1}^{n}\left\Vert
u\right\Vert _{H^{2,1}\left( \Psi _{iT}^{\pm }\right) }^{2}<\infty \right\} ,
\\ 
H^{1,0}\left( S_{T}\right) =\left\{ 
\begin{array}{c}
u:\left\Vert u\right\Vert _{H^{1,0}\left( S_{T}\right) }^{2}= \\ 
=\sum\limits_{i=1}^{n}\left( \sum\limits_{j=1,j\neq i}^{n}\left\Vert
u_{x_{j}}\right\Vert _{L_{2}\left( \Psi _{iT}^{\pm }\right) }^{2}+\left\Vert
u\right\Vert _{L_{2}\left( \Psi _{iT}^{\pm }\right) }^{2}\right) <\infty%
\end{array}%
\right\} .%
\end{array}%
\right.
\end{equation*}

The Mean Field Games System (MFGS) is \cite{A}:%
\begin{equation}
\left. 
\begin{array}{c}
u_{t}(x,t)+\Delta u(x,t){-k(x)(\nabla u(x,t))^{2}/2}+ \\ 
+\int\limits_{\Omega }K\left( x,y\right) m\left( y,t\right) dy+f\left(
x,t\right) m\left( x,t\right) =0,\text{ }\left( x,t\right) \in Q_{T}, \\ 
m_{t}(x,t)-\Delta m(x,t){-\text{div}(k(x)m(x,t)\nabla u(x,t))}=0,\text{ }%
\left( x,t\right) \in Q_{T}.%
\end{array}%
\right.  \label{2.2}
\end{equation}%
Here $u(x,t)$ is the value function, and $m(x,t)$ is the density of
agents/players at the point $x$ and at the moment of time $t$. The meaning
of the coefficient ${k(x)}$ was described in Introduction. The integral term
in (\ref{2.2}) is the so-called \textquotedblleft global interaction" term.
This term expresses the action on an agent occupying the point $x$ by the
rest of agents \cite[page 634]{MFG1}. The first and second equations (\ref%
{2.2}) are Hamilton-Jacobi-Bellman (HJB) and Fokker-Planck (FP) equations
respectively.

When working with a forward problem for system (\ref{2.2}), one usually
imposes some boundary conditions, initial condition $m\left( x,0\right) $
and terminal condition $u\left( x,T\right) $ \cite{A,LL2}. However, since we
are interested in an inverse problem, we do not discuss the forward problem
here. Rather we formulate now the CIP, which we study in this paper. Let $%
t_{0}\in \left( 0,T\right) $\ be a number. For brevity of further
developments, we set $t_{0}=T/2.$

\textbf{Coefficient Inverse} \textbf{Problem }(CIP). \emph{Assume that
functions }$u,m\in C^{4}\left( \overline{Q}_{T}\right) $\emph{\ satisfy
equations (\ref{2.2}). Let}%
\begin{equation}
\left. 
\begin{array}{c}
u\left( x,T/2\right) =u_{0}\left( x\right) ,\text{ }m\left( x,T/2\right)
=m_{0}\left( x\right) ,\text{ }x\in \Omega , \\ 
u\mid _{S_{T}}=g_{0}\left( x,t\right) ,\text{ }u_{x_{1}}\mid _{\Psi
_{1T}^{+}}=g_{1}\left( x,t\right) , \\ 
m\mid _{S_{T}}=p_{0}\left( x,t\right) ,\text{ }m_{x_{1}}\mid _{\Psi
_{1,T}^{+}}=p_{1}\left( x,t\right) ,%
\end{array}%
\right.  \label{2.3}
\end{equation}%
\emph{Determine the coefficient }$k\left( x\right) \in C^{1}\left( \overline{%
\Omega }\right) ,$\emph{\ assuming that functions }$g_{0},g_{1},p_{0},p_{1}$%
\emph{\ in (\ref{2.3}) are known. }

\textbf{Remarks 2.1:}

\emph{1.\ \ As to the required smoothness }$u,m\in C^{4}\left( \overline{Q}%
_{T}\right) ,$\emph{\ we note that the minimal smoothness conditions are
traditionally a minor concern in the theory of CIPs, see, e.g. \cite{Nov,Rom}%
.}

\emph{2. Note that while the Dirichlet data }$g_{0}$\emph{\ and }$p_{0}$%
\emph{\ in (\ref{2.3}) are given on the entire lateral boundary }$S_{T},$%
\emph{\ the Neumann data are given only on the part }$\Psi _{1,T}^{+}$\emph{%
\ of this boundary. In other words, we work with partially incomplete data.}

\subsection{Discussion}

\label{sec:2.2}

CIP\ (\ref{2.3}), (\ref{2.3}) is the CIP\ with the single measurement data.
In an applied case, the input data (\ref{2.3}) can be obtained by conducting
polls of game players. In the case of the boundary data, polls are usually
conducted not just at the boundary but rather in a small neighborhood of the
boundary. This way the normal derivatives $u_{x}$ and $m_{x}$ at the part of
the boundary $\Psi _{1,T}^{+}$ are approximately computed. To find functions 
$u\left( x,T/2\right) $ and $m\left( x,T/2\right) $ in (\ref{2.3}), one
needs to have polling inside of the domain $\Omega $ at the moment of time $%
t=T/2$ ($t=t_{0}$ actually: recall that we set $t_{0}=T/2$ for brevity).
These are not excessive requirements for the volume of the input data. To
compare, we mention that it is assumed in \cite{Ding} that the solution of
the MFGS is known for all $\left( x,t\right) \in Q_{T}.$ Our method provides
error estimates of the solution with respect to the noise in the data. In
addition, H\"{o}lder stability estimate for this CIP was obtained in \cite[%
Theorem 3.3]{MFG6}, and they imply uniqueness of this problem.

The main challenges of working with the MFGS (\ref{2.2}) are:

1. The presence of the integral operator in the HJB equation. Similar terms
are not present in any past publications for CIPs for a single parabolic
equation.

2. The presence of the term $\Delta u$ in the FP equation. This
significantly complicates the analysis of system (\ref{2.2}) since $\Delta u$
is involved in the principal part $u_{t}+\Delta u$ of the PDE operator of
the HJB equation.

3. The nonlinearity of both PDEs (\ref{2.2}).

4. The fact that times are running in two opposite directions in these
equations, which means that the classical theory of parabolic equations \cite%
{Lad} is inapplicable here.

The apparatus of Carleman estimates allows us to handle these challenges.

\subsection{The form of the kernel $K\left( x,y\right) $}

\label{sec:2.3}

We will work with the same two forms of the kernel $K\left( x,y\right) $ as
the ones in \cite{MFG6}. It was proposed in \cite[formula (2.7)]{Ding} to
choose $K\left( x,y\right) $ as the product of Gaussians. Since the $\delta
- $ function is the limiting case of the Gaussian, then the first form of $%
K\left( x,y\right) ,$ which we use in this paper is:%
\begin{equation}
\left. 
\begin{array}{c}
K_{1}\left( x,y\right) =\delta \left( x_{1}-y_{1}\right) Y_{1}\left( x,%
\overline{y}\right) , \\ 
Y_{1}\left( x,\overline{y}\right) \in L_{\infty }\left( \Omega \times \Omega
_{1}\right) ,\mbox{ }\left\Vert Y_{1}\right\Vert _{L_{\infty }\left( \Omega
\times \Omega _{1}\right) }\leq M,%
\end{array}%
\right.  \label{2.4}
\end{equation}%
where $M>0$ is a number. Let $H\left( s\right) $, $s\in \mathbb{R}$ be the
Heaviside function, 
\begin{equation*}
H\left( s\right) =\left\{ 
\begin{array}{c}
1,\text{ }s>0, \\ 
0,\text{ }s<0.%
\end{array}%
\right.
\end{equation*}%
The second form of $K\left( x,y\right) $ used here is:%
\begin{equation}
\left. 
\begin{array}{c}
K_{2}\left( x,y\right) =H\left( y_{1}-x_{1}\right) Y_{2}\left( x,y\right) ,%
\text{ }\left( x,y\right) \in \Omega \times \Omega , \\ 
\text{ }Y_{2}\in L_{\infty }\left( \Omega \times \Omega _{1}\right) ,\text{ }%
\left\Vert \text{ }Y_{2}\right\Vert _{L_{\infty }\left( \Omega \times \Omega
_{1}\right) }\leq M.%
\end{array}%
\right.  \label{2.5}
\end{equation}%
Using (\ref{2.1}), we obtain that in the case (\ref{2.4})%
\begin{equation}
\int\limits_{\Omega }K_{1}\left( x,y\right) m\left( y,t\right)
dy=\int\limits_{\Omega _{1}}Y_{1}\left( x,\overline{y}\right) m\left( x_{1},%
\overline{y},t\right) d\overline{y}dt,\text{ }\left( x,t\right) \in Q_{T}.
\label{2.6}
\end{equation}%
Similarly, in the case (\ref{2.5}) 
\begin{equation}
\left. 
\begin{array}{c}
\int\limits_{\Omega }K_{2}\left( x,y\right) m\left( y,t\right) dy= \\ 
=\int\limits_{\Omega _{1}}\left( \int\limits_{x_{1}}^{b}Y_{2}\left( x_{1},%
\overline{x},y_{1},\overline{y}\right) m\left( y_{1},\overline{y},t\right)
dy_{1}\right) d\overline{y},\text{ }\left( x,t\right) \in Q_{T}.%
\end{array}%
\right.  \label{2.7}
\end{equation}

\section{The Transformation Procedure}

\label{sec:3}

This procedure transforms CIP (\ref{2.2}), (\ref{2.3}) to a BVP for a system
of four integral differential equations, which does not contain the unknown
coefficient $k\left( x\right) .$ On the other hand, the transformed system
contains Volterra integrals with respect to $t$.

Recall that the function $u_{0}\left( x\right) =u\left( x,T/2\right) $ is
given (\ref{2.3}) as a representative of the given data in (\ref{2.3}). We
assume that 
\begin{equation}
\left\vert \nabla u_{0}\left( x\right) \right\vert ^{2}\geq c>0,\text{ in }%
\overline{\Omega },  \label{3.1}
\end{equation}%
where $c$ is a number. Set in the first equation (\ref{2.2}) $t=T/2.$ Using (%
\ref{2.3}),\emph{\ }we obtain 
\begin{equation*}
k\left( x\right) =\frac{2\left[ u_{t}\left( x,T/2\right) +\Delta u_{0}\left(
x\right) +\int\limits_{\Omega }K\left( x,y\right) m_{0}\left( y\right)
dy+f\left( x,T/2\right) m_{0}\left( x\right) \right] }{\left\vert \nabla
u_{0}\left( x\right) \right\vert ^{2}}.
\end{equation*}%
This is equivalent with%
\begin{align}
k\left( x\right) &=\frac{2}{\left\vert \nabla u_{0}\left( x\right)
\right\vert ^{2}}\left[ u_{t}\left( x,T/2\right) +F\left( x\right) \right] ,
\label{3.2} \\
F\left( x\right) &=\Delta u_{0}\left( x\right) +\int\limits_{\Omega }K\left(
x,y\right) m_{0}\left( y\right) dy+f\left( x,T/2\right) m_{0}\left( x\right)
.  \label{3.3}
\end{align}%
Substituting (\ref{3.2}) and (\ref{3.3}) in (\ref{2.2}), we obtain two
equations in which the unknown coefficient $k\left( x\right) $ is not
involved. The first equation is:

\begin{equation}
\left. 
\begin{array}{c}
u_{t}(x,t)+\Delta u(x,t){-(\nabla u(x,t))^{2}}\left[ u_{t}\left(
x,T/2\right) +F\left( x\right) \right] \left\vert \nabla u_{0}\left(
x\right) \right\vert ^{-2}+ \\ 
+\int\limits_{\Omega }K\left( x,y\right) m\left( y,t\right) dy+f\left(
x,t\right) m\left( x,t\right) =0,\text{ }\left( x,t\right) \in Q_{T}.%
\end{array}%
\right.  \label{3.4}
\end{equation}%
The second equation is:%
\begin{equation}
\left. 
\begin{array}{c}
m_{t}(x,t)-\Delta m(x,t){-} \\ 
-{2\text{div}}\left\{ \left[ u_{t}\left( x,T/2\right) +F\left( x\right) %
\right] \left\vert \nabla u_{0}\left( x\right) \right\vert ^{-2}{%
m(x,t)\nabla u(x,t)}\right\} =0,\text{ }\left( x,t\right) \in Q_{T}.%
\end{array}%
\right.  \label{3.5}
\end{equation}

Even though equations (\ref{3.4}) and (\ref{3.5}) do not contain the unknown
coefficient $k\left( x\right) $, still the presence of the function $%
u_{t}\left( x,T/2\right) $ in them makes these equations non-local ones.
Thus, we proceed further. Differentiate both parts of these equations with
respect to $t$ and denote:%
\begin{align}
v\left( x,t\right) &=u_{t}\left( x,t\right) ,\text{ }p\left( x,t\right)
=m_{t}\left( x,t\right) ,\text{ }\left( x,t\right) \in Q_{T},  \label{3.6} \\
w\left( x,t\right) &=v_{t}\left( x,t\right) ,\text{ }q\left( x,t\right)
=p_{t}\left( x,t\right) ,\text{ }\left( x,t\right) \in Q_{T}.  \label{3.7}
\end{align}%
By (\ref{2.3}) and (\ref{3.6}) 
\begin{equation}
\left. 
\begin{array}{c}
u\left( x,t\right) =\int\limits_{T/2}^{t}v\left( x,\tau \right) d\tau
+u_{0}\left( x\right) ,\text{ } \\ 
m\left( x,t\right) =\int\limits_{T/2}^{t}p\left( x,\tau \right) d\tau
+m_{0}\left( x\right) ,\text{ }\left( x,t\right) \in Q_{T}.%
\end{array}%
\right.  \label{3.8}
\end{equation}%
Next, by (\ref{3.6}) and (\ref{3.7}) 
\begin{equation}
u_{t}\left( x,T/2\right) =v\left( x,T/2\right) =v\left( x,t\right)
-\int\limits_{T/2}^{t}w\left( x,\tau \right) d\tau .  \label{3.9}
\end{equation}

Hence, using (\ref{3.4})-(\ref{3.9}), we obtain two integral differential
equations. The first equation is:%
\begin{equation}
\left. 
\begin{array}{c}
L_{1}\left( v,w,p,q\right) = \\ 
=v_{t}(x,t)+\Delta v(x,t){-2\nabla v(x,t)\cdot }\left(
\int\limits_{T/2}^{t}\nabla v\left( x,\tau \right) d\tau +\nabla u_{0}\left(
x\right) \right) \times \\ 
\times \left[ \left( v\left( x,t\right) -\int\limits_{T/2}^{t}w\left( x,\tau
\right) d\tau \right) +F\left( x\right) \right] \left\vert \nabla
u_{0}\left( x\right) \right\vert ^{-2}+ \\ 
+\int\limits_{\Omega }K\left( x,y\right) p\left( y,t\right) dy+f\left(
x,t\right) p\left( x,t\right) + \\ 
+f_{t}\left( x,t\right) \left( \int\limits_{T/2}^{t}p\left( x,\tau \right)
d\tau +m_{0}\left( x\right) \right) =0,\text{ }\left( x,t\right) \in Q_{T}.%
\end{array}%
\right.  \label{3.10}
\end{equation}%
The second equation is:%
\begin{equation}
\left. 
\begin{array}{c}
L_{2}\left( v,w,p,q\right) = \\ 
=p_{t}(x,t)-\Delta p(x,t){-} \\ 
-{2\text{div}}\left\{ 
\begin{array}{c}
\left[ v\left( x,t\right) -\int\limits_{T/2}^{t}w\left( x,\tau \right) d\tau
+F\left( x\right) \right] \left\vert \nabla u_{0}\left( x\right) \right\vert
^{-2}\times \\ 
\times p{(x,t)}\left( \int\limits_{T/2}^{t}\nabla v\left( x,\tau \right)
d\tau +\nabla u_{0}\left( x\right) \right)%
\end{array}%
\right\} - \\ 
-2\text{div}\left\{ 
\begin{array}{c}
\left[ v\left( x,t\right) -\int\limits_{T/2}^{t}w\left( x,\tau \right) d\tau
+F\left( x\right) \right] \left\vert \nabla u_{0}\left( x\right) \right\vert
^{-2}\times \\ 
\times {\nabla v(x,t)}\left( \int\limits_{T/2}^{t}p\left( x,\tau \right)
d\tau +m_{0}\left( x\right) \right)%
\end{array}%
\right\} =0,\text{ } \\ 
\left( x,t\right) \in Q_{T}.%
\end{array}%
\right.  \label{3.11}
\end{equation}%
An inconvenient point of equation (\ref{3.11}), which does not allow to
apply a Carleman estimate, is the presence of the terms with the mixed $x,t$
derivatives of the function $v$, 
\begin{equation*}
\int\limits_{T/2}^{t}\partial _{x_{i}}w\left( x,\tau \right) d\tau
=\int\limits_{T/2}^{t}\partial _{x_{i}}v_{t}\left( x,\tau \right) d\tau
\end{equation*}

Hence, we differentiate both equations (\ref{3.10}) and (\ref{3.11}) with
respect to $t$ and use (\ref{3.6}) and (\ref{3.7}). Then we obtain from (\ref%
{3.10}): 
\begin{equation}
\left. 
\begin{array}{c}
L_{3}\left( v,w,p,q\right) = \\ 
=w_{t}(x,t)+\Delta w(x,t){-} \\ 
-2\left[ {\nabla w(x,t)\cdot }\left( \int\limits_{T/2}^{t}\nabla v\left(
x,\tau \right) d\tau +\nabla u_{0}\left( x\right) \right) +\left( \nabla
v\right) ^{2}{(x,t)}\right] \times \\ 
\times \left[ \left( v\left( x,t\right) -\int\limits_{T/2}^{t}w\left( x,\tau
\right) d\tau \right) +F\left( x\right) \right] \left\vert \nabla
u_{0}\left( x\right) \right\vert ^{-2}+ \\ 
+\int\limits_{\Omega }K\left( x,y\right) q\left( y,t\right) dy+f\left(
x,t\right) q\left( x,t\right) +2f_{t}\left( x,t\right) p\left( x,t\right) +
\\ 
+f_{tt}\left( x,t\right) \left( \int\limits_{T/2}^{t}p\left( x,\tau \right)
d\tau +m_{0}\left( x\right) \right) =0,\text{ }\left( x,t\right) \in Q_{T}.%
\end{array}%
\right.  \label{3.12}
\end{equation}%
And we obtain from (\ref{3.11}):%
\begin{equation}
\left. 
\begin{array}{c}
L_{4}\left( v,w,p,q\right) = \\ 
=q_{t}(x,t)-\Delta q(x,t){-} \\ 
-{2\text{div}}\left\{ 
\begin{array}{c}
\left[ v\left( x,t\right) -\int\limits_{T/2}^{t}w\left( x,\tau \right) d\tau
+F\left( x\right) \right] \left\vert \nabla u_{0}\left( x\right) \right\vert
^{-2}\times \\ 
\times q{(x,t)}\left( \int\limits_{T/2}^{t}\nabla v\left( x,\tau \right)
d\tau +\nabla u_{0}\left( x\right) \right) +p{(x,t)\nabla v}%
\end{array}%
\right\} - \\ 
-2\text{div}\left\{ 
\begin{array}{c}
\left[ v\left( x,t\right) -\int\limits_{T/2}^{t}w\left( x,\tau \right) d\tau
+F\left( x\right) \right] \left\vert \nabla u_{0}\left( x\right) \right\vert
^{-2}\times \\ 
\times \left[ 
\begin{array}{c}
{\nabla w(x,t)}\left( \int\limits_{T/2}^{t}p\left( x,\tau \right) d\tau
+m_{0}\left( x\right) \right) \\ 
+{\nabla v(x,t)p(x,t)}%
\end{array}%
\right]%
\end{array}%
\right\} =\text{ } \\ 
=0,\text{ }\left( x,t\right) \in Q_{T}.%
\end{array}%
\right.  \label{3.13}
\end{equation}

Boundary conditions for the vector function $\left( v,w,p,q\right) $ are
derived from (\ref{2.3}), (\ref{3.6}) and (\ref{3.7}), 
\begin{equation}
\left. 
\begin{array}{c}
v\mid _{S_{T}}=\partial _{t}g_{0}\left( x,t\right) ,\text{ }v_{x_{1}}\mid
_{\Psi _{1T}^{+}}=\partial _{t}g_{1}\left( x,t\right) , \\ 
w\mid _{S_{T}}=\partial _{t}^{2}g_{0}\left( x,t\right) ,\text{ }%
w_{x_{1}}\mid _{\Psi _{1T}^{+}}=\partial _{t}^{2}g_{1}\left( x,t\right) , \\ 
p\mid _{S_{T}}=\partial _{t}p_{0}\left( x,t\right) ,\text{ }p_{x_{1}}\mid
_{\Psi _{1T}^{+}}=\partial _{t}p_{1}\left( x,t\right) , \\ 
q\mid _{S_{T}}=\partial _{t}^{2}p_{0}\left( x,t\right) ,\text{ }%
q_{x_{1}}\mid _{\Psi _{1T}^{+}}=\partial _{t}^{2}p_{1}\left( x,t\right) .%
\end{array}%
\right.  \label{3.14}
\end{equation}

This concludes the transformation procedure. Thus, we need to solve below
BVP (\ref{3.10})-(\ref{3.14}). If this problem is solved, then the unknown
coefficient $k\left( x\right) $ should be computed by a slightly modified
formula (\ref{3.2}), where $F\left( x\right) $ is given in (\ref{3.3}) 
\begin{equation}
k\left( x\right) =\frac{2}{\left\vert \nabla u_{0}\left( x\right)
\right\vert ^{2}}\left[ v\left( x,T/2\right) +F\left( x\right) \right] .
\label{3.15}
\end{equation}

\section{Convexification}

\label{sec:4}

\subsection{The Carleman Weight Function (CWF) and some estimates}

\label{sec:4.1}

First, we introduce our CWF. Since the variable $x_{1}$ is singled out in (%
\ref{2.4})-(\ref{2.7}), then our CWF depends only on one spatial variable,
which is $x_{1}.$ Our CWF is:%
\begin{equation}
\varphi _{\lambda }\left( x,t\right) =\varphi _{\lambda }\left(
x_{1},t\right) =\exp \left[ 2\lambda \left( x_{1}^{2}-\left( t-T/2\right)
^{2}\right) \right] ,  \label{4.1}
\end{equation}%
where $\lambda \geq 1$ is a large parameter. We will choose $\lambda $
later. By (\ref{2.1}) and (\ref{4.1})%
\begin{equation}
\left. 
\begin{array}{c}
\max_{\overline{Q}_{T}}\varphi _{\lambda }\left( x_{1},t\right) =e^{2\lambda
b^{2}},\text{ } \\ 
\max_{x_{1}\in \left[ a,b\right] }\varphi _{\lambda }\left( x_{1},T\right)
=\max_{x_{1}\in \left[ a,b\right] }\varphi _{\lambda }\left( x_{1},0\right)
=e^{2\lambda b^{2}}e^{-\lambda T^{2}/2}, \\ 
\min_{\overline{Q}_{\gamma ,T}}\varphi _{\lambda }\left( x_{1},t\right)
=\exp \left[ 2\lambda \left( a^{2}-\left( \gamma T\right) ^{2}/4\right) %
\right] \geq \exp \left( -\lambda \left( \gamma T\right) ^{2}/2\right) .%
\end{array}%
\right.  \label{4.2}
\end{equation}

\textbf{Lemma 4.1} \cite{MFG6}. \emph{The following estimates of integrals (%
\ref{2.6}) and (\ref{2.7}) are valid: }%
\begin{equation*}
\left. 
\begin{array}{c}
\int\limits_{Q_{T}}\left( \int\limits_{\Omega _{1}}Y_{1}\left( x,\overline{y}%
\right) s\left( x_{1},\overline{y},t\right) d\overline{y}\right) ^{2}\varphi
_{\lambda }\left( x_{1},t\right) dxdt\leq \\ 
\leq C_{1}\int\limits_{Q_{T}}s^{2}\varphi _{\lambda }dxdt,\text{ }\forall
s\in L_{2}\left( Q_{T}\right) ,\text{ }\forall \lambda >0.%
\end{array}%
\right.
\end{equation*}%
\emph{And also}%
\begin{equation*}
\left. 
\begin{array}{c}
\int\limits_{Q_{T}}\left[ \int\limits_{\Omega _{1}}\left(
\int\limits_{x_{1}}^{b}Y_{2}\left( x_{1},\overline{x},y_{1},\overline{y}%
\right) s\left( y_{1},\overline{y},t\right) dy_{1}\right) d\overline{y}%
\right] ^{2}\varphi _{\lambda }\left( x_{1},t\right) dxdt\leq \\ 
\leq C_{1}\int\limits_{Q_{T}}s^{2}\varphi _{\lambda }dxdt,\text{ }\forall
s\in L_{2}\left( Q_{T}\right) ,\text{ }\forall \lambda >0,%
\end{array}%
\right.
\end{equation*}%
\emph{where the number }$C_{1}=C_{1}\left( \Omega ,M\right) >0$\emph{\
depends only on }$\Omega $\emph{\ and }$M.$

\textbf{Lemma 4.2 }\cite{Ksurvey}, \textbf{\ }\cite[Lemma 3.1.1]{KL}.\textbf{%
\ }\emph{The following estimate is valid:}%
\begin{equation*}
\int\limits_{Q_{T}}\left( \int\limits_{T/2}^{t}s\left( x,\tau \right) d\tau
\right) ^{2}\varphi _{\lambda }dxdt\leq \frac{C_{2}}{\lambda }%
\int\limits_{Q_{T}}\left( s^{2}\varphi _{\lambda }\right) \left( x,t\right)
dxdt,\text{ }\forall \lambda >0,\text{ }\forall s\in L_{2}\left(
Q_{T}\right) ,
\end{equation*}%
\emph{\ \ where the number }$C_{2}>0$\emph{\ is independent from }$\lambda
,s.$

\textbf{Lemma 4.3} (Two Carleman estimates \cite{MFG6}) \emph{Let the domain 
}$\Omega $\emph{\ be as in (\ref{2.1}) and }$\varphi _{\lambda }\left(
x,t\right) $\emph{\ be the function defined in (\ref{4.1}). Then there
exists a sufficiently large number }$\overline{\lambda }_{0}=\overline{%
\lambda }_{0}\left( \Omega ,T\right) \geq 1$\emph{\ and a number }$%
C_{3}=C_{3}\left( \Omega ,T\right) >0,$\emph{\ both numbers depend only on }$%
\Omega $ \emph{and} $T$, \emph{such that the following two Carleman
estimates are valid: }%
\begin{equation}
\left. 
\begin{array}{c}
\int\limits_{Q_{T}}\left( u_{t}\pm \Delta u\right) ^{2}\varphi _{\lambda
}dxdt\geq \left( C_{3}/\lambda \right) \int\limits_{Q_{T}}\left(
u_{t}^{2}+\sum\limits_{i,j=1}^{n}u_{x_{i}x_{j}}^{2}\right) \varphi _{\lambda
}dxdt+ \\ 
+C_{3}\lambda \int\limits_{Q_{T}}\left( \left( \nabla u\right) ^{2}+\lambda
^{2}u^{2}\right) \varphi _{\lambda }dxdt- \\ 
-C_{3}\left( \left\Vert u\left( x,T\right) \right\Vert _{H^{1}\left( \Omega
\right) }^{2}+\left\Vert u\left( x,0\right) \right\Vert _{H^{1}\left( \Omega
\right) }^{2}\right) \exp \left( -\lambda T^{2}/2+2\lambda b^{2}\right) ,%
\text{ } \\ 
\forall \lambda \geq \overline{\lambda },\text{ }\forall u\in H^{2,1}\left(
Q_{T}\right) \cap \left\{ u:u\mid _{S_{T}}=u_{x_{1}}\mid _{\Psi
_{1,T}^{+}{}}=0\right\} .%
\end{array}%
\right.  \label{4.3}
\end{equation}

\subsection{The functional}

To solve BVP (\ref{3.10})-(\ref{3.14}), we now construct the convexified
cost functional. First, we introduce the number $k_{n},$ 
\begin{equation}
k_{n}=\left[ \frac{n+1}{2}\right] +2,  \label{4.04}
\end{equation}%
where $\left[ \left( n+1\right) /2\right] $ is the maximal integer, which
does not exceed $\left( n+1\right) /2.$ By embedding theorem, $%
H^{k_{n}}\left( Q_{T}\right) \subset C^{1}\left( \overline{Q}_{T}\right) $
as a set and%
\begin{equation}
\left\Vert z\right\Vert _{C^{1}\left( \overline{Q}_{T}\right) }\leq
C\left\Vert z\right\Vert _{H^{k_{n}}\left( Q_{T}\right) },\text{ }\forall
z\in H^{k_{n}}\left( Q_{T}\right) ,  \label{4.4}
\end{equation}%
where the number $C=C\left( Q_{T}\right) $ depends only on the domain $%
Q_{T}, $ also, see the first item of Remarks 2.1. Introduce two Hilbert
spaces:%
\begin{equation}
\left. 
\begin{array}{c}
H=\left( H^{k_{n}}\left( Q_{T}\right) \right) ^{4}= \\ 
=\left\{ U=\left( u_{1},u_{2},u_{3},u_{4}\right) :\left\Vert U\right\Vert
_{H}^{2}=\sum\limits_{i=1}^{4}\left\Vert u_{i}\right\Vert _{H^{k_{n}}\left(
Q_{T}\right) }^{2}<\infty \right\} , \\ 
H_{0}=\left\{ 
\begin{array}{c}
U=\left( u_{1},u_{2},u_{3},u_{4}\right) \in H: \\ 
u_{i}\mid _{S_{T}}=\partial _{x_{1}}u_{i}\mid _{\Psi _{1,T}^{+}}=0,i=1,...,4%
\end{array}%
\right\} .%
\end{array}%
\right.  \label{5.1}
\end{equation}%
Below $\left[ ,\right] $ is the scalar product in $H$.

Let $R>0$ be an arbitrary number. Denote $U=\left( v,w,p,q\right) .$ Define
the set $B\left( R\right) ,$%
\begin{equation}
\left. B\left( R\right) =\left\{ 
\begin{array}{c}
U=\left( v,w,p,q\right) \in H:\left\Vert U\right\Vert _{H}<R, \\ 
\text{functions }v,w,p,q\text{ satisfy boundary conditions (\ref{3.14})}%
\end{array}%
\right\} .\right.  \label{4.5}
\end{equation}%
We solve problem (\ref{3.10})-(\ref{3.14}) via the solution of the following
problem:

\textbf{Minimization Problem:}\emph{\ Minimize the functional }$J_{\lambda
,\beta }\left( U\right) $\emph{\ on the set }$\overline{B\left( R\right) },$%
\emph{\ where}%
\begin{equation}
\left. 
\begin{array}{c}
J_{\lambda ,\beta }\left( U\right) = \\ 
=e^{-2\lambda b^{2}}\int\limits_{Q_{T}}\lambda ^{3/2}\left[ \left(
L_{1}\left( U\right) \right) ^{2}+\left( L_{3}\left( U\right) \right) ^{2}%
\right] \varphi _{\lambda }\left( x,t\right) dxdt+ \\ 
+e^{-2\lambda b^{2}}\int\limits_{Q_{T}}\left[ \left( L_{2}\left( U\right)
\right) ^{2}+\left( L_{4}\left( U\right) \right) ^{2}\right] \varphi
_{\lambda }\left( x,t\right) dxdt+\beta \left\Vert U\right\Vert _{H}^{2}.%
\end{array}%
\right.  \label{4.6}
\end{equation}%
\emph{\ }

In (\ref{4.6}), $\beta \in \left( 0,1\right) $ is the regularization
parameter and $e^{-2\lambda b^{2}}$ is the balancing factor between integral
terms and the regularization term. Indeed, by (\ref{4.2}) $e^{-2\lambda
b^{2}}\varphi _{\lambda }\left( x,t\right) \leq 1$ in $\overline{Q}_{T},$
which balances the integral terms with $\beta \in \left( 0,1\right) .$

\section{Convergence Analysis}

\label{sec:5}

\subsection{The central theorem}

\label{sec:5.1}

\textbf{Theorem 5.1} (the central result). \emph{Assume that in (\ref{2.3})
functions }$u_{0},m_{0}\in H^{2}\left( \Omega \right) $, 
\begin{equation}
\emph{\ \ }\left\Vert u_{0}\right\Vert _{H^{2}\left( \Omega \right) },%
\mbox{
	}\left\Vert m_{0}\right\Vert _{H^{2}\left( \Omega \right) }<R  \label{5.02}
\end{equation}%
\emph{and condition (\ref{3.1}) holds. In addition, let either of two
conditions (\ref{2.4}) or (\ref{2.5}) holds. Let }$B\left( R\right) $ \emph{%
be the set defined in (\ref{4.5}). Then:}

\emph{1. The functional }$J_{\lambda ,\beta }\left( U\right) $\emph{\ has Fr%
\'{e}chet derivative }$J_{\lambda ,\beta }^{\prime }\left( U\right) \in
H_{0} $\emph{\ at every point }$U\in \overline{B\left( R\right) },$\emph{\
and this derivative is Lipschitz continuous on }$\overline{B\left( R\right) }
$\emph{, i.e. there exists a number }$D=D\left( \Omega ,T,M,c,R,\lambda
\right) >0$\emph{\ depending only on listed parameters such that \ }%
\begin{equation}
\left\Vert J_{\lambda ,\beta }^{\prime }\left( U_{2}\right) -J_{\lambda
,\beta }^{\prime }\left( U_{1}\right) \right\Vert _{H}\leq D\left\Vert
U_{2}-U_{1}\right\Vert _{H},\text{ }\forall U_{1},U_{2}\in \overline{B\left(
R\right) }.  \label{5.3}
\end{equation}%
\emph{2. Let }$\overline{\lambda }_{0}=\overline{\lambda }_{0}\left( \Omega
,T\right) \geq 1$\emph{\ be the number of Lemma 4.3. There exist a
sufficiently large number }$\overline{\lambda }_{1}=\overline{\lambda }%
_{1}\left( \Omega ,T,M,c,R\right) \geq \overline{\lambda }_{0}$\emph{\ and a
number }$C_{4}=C_{4}\left( \Omega ,T,M,c,R\right) >0,$\emph{\ both numbers
depend only on listed parameters, such that if }$\lambda \geq \overline{%
\lambda }_{1}$\emph{\ and }$\beta \in \left[ 2\exp \left( -\lambda
T^{2}/4\right) ,1\right) ,$\emph{\ then the functional }$J_{\lambda ,\beta }$%
\emph{\ is strongly convex on }$\overline{B\left( R\right) },$\emph{\ i.e.
the following inequality holds:}%
\begin{equation}
\left. 
\begin{array}{c}
J_{\lambda ,\beta }\left( U_{2}\right) -J_{\lambda ,\beta }\left(
U_{1}\right) -\left[ J_{\lambda ,\beta }^{\prime }\left( U_{1}\right)
,U_{2}-U_{1}\right] \geq \\ 
\geq C_{4}\exp \left[ -\lambda \left( \gamma T\right) ^{2}/2\right]
\left\Vert U_{2}-U_{1}\right\Vert _{\left( H^{2,1}\left( Q_{\gamma T}\right)
\right) ^{4}}^{2}+\left( \beta /2\right) \left\Vert U_{2}-U_{1}\right\Vert
_{H}^{2}, \\ 
U_{1},U_{2}\in \overline{B\left( R\right) },\forall \geq \overline{\lambda }%
_{1},\forall \gamma \in \left( 0,1\right) .%
\end{array}%
\right.  \label{5.4}
\end{equation}%
\emph{3. For }$\lambda $\emph{\ and }$\beta $\emph{\ as in item 2, the
functional }$J_{\lambda ,\beta }$\emph{\ has unique minimizer }$U_{\min
,\lambda ,\beta }\in $\emph{\ }$\overline{B\left( R\right) }$ \emph{on this
set,\ and the following inequality holds}$:$\emph{\ }%
\begin{equation}
\hspace{-1.5cm}\left[ J_{\lambda ,\beta }^{\prime }\left( U_{\min ,\lambda
,\beta }\right) ,\left( U_{\min ,\lambda ,\beta }-U\right) \right] \leq 0,%
\mbox{ }\forall U\in \overline{B\left( R\right) }.  \label{5.5}
\end{equation}

\textbf{Remarks 5.1:}

\begin{enumerate}
\item \emph{Below }$C_{4}>0$\emph{\ denotes different numbers depending only
on parameters indicated above.}

\item \emph{Theorem 5.1 works only for sufficiently large values of the
parameter }$\lambda $\emph{. However, it was established in our past works
on the convexification that it actually works well for }$\lambda \in \left[
1,5\right] $\emph{\ \cite{Bak,Kpar,KL,MFG8}. In particular, in the current
paper }$\lambda =3$\emph{\ provides accurate solutions, see Test 1 in
section 7 for the choice of an optimal value of }$\lambda $\emph{. Basically
the same thing takes place in any asymptotic theory. Indeed, such a theory
usually states that if a certain parameter }$X$\emph{\ is sufficiently
large, then a formula }$Y$\emph{\ provides a good accuracy. However, in
computations only a numerical experience can establish reasonable values of }%
$X.$
\end{enumerate}

\subsection{Accuracy estimate and global convergence of the gradient descent
method}

\label{sec:5.2}

The next two natural questions after Theorem 5.1 are about an estimate of
the accuracy of the minimizer of the functional $J_{\lambda ,\beta }$ on the
set $\overline{B\left( R\right) }$, which was found in that theorem, and
also about the global convergence of the gradient descent method of the
minimization of $J_{\lambda ,\beta }.$ We address these two questions in
this subsection. This is done quite similarly with \cite{MFG8}. Hence, we
omit proofs of Theorems 5.2 and 5.3 since they are very similar with the
proofs of Theorems 4.4 and 4.5 of \cite{MFG8} respectively.

\subsubsection{Accuracy estimate}

Suppose that our boundary data (\ref{3.14}) are given with a noise of the
level $\delta \in \left( 0,1\right) $. By one of principles of the theory of
Ill-Posed Problems \cite{T}, we assume that there exists an exact solution 
\begin{equation}
\left( k^{\ast },v^{\ast },w^{\ast },p^{\ast },q^{\ast }\right) =\left(
k^{\ast },U^{\ast }\right) \in C^{1}\left( \overline{\Omega }\right) \times
B^{\ast }\left( R-\delta \right)  \label{5.50}
\end{equation}%
of our CIP (\ref{2.2}), (\ref{2.3}) with the exact, i.e. noiseless data.
Uniqueness of this solution follows from the uniqueness theorem for our CIP\ 
\cite{MFG6}. The exact data are equipped below by the superscript
\textquotedblleft $^{\ast }$". So, the set $B^{\ast }\left( R\right) $ is
defined similarly with the set $B\left( R\right) $ in (\ref{4.5})%
\begin{equation}
\left. B^{\ast }\left( R\right) =\left\{ 
\begin{array}{c}
U=\left( v,w,p,q\right) \in H:\left\Vert U\right\Vert _{H}<R, \\ 
\text{functions }v,w,p,q\text{ satisfy boundary conditions (\ref{3.14}) } \\ 
\text{with the exact data }\partial ^{k}g_{0}^{\ast },\partial
^{k}g_{1}^{\ast },\partial ^{k}p_{0}^{\ast },\partial ^{k}p_{1}^{\ast
},k=1,2.%
\end{array}%
\right\} .\right.  \label{5.6}
\end{equation}%
Assume that there exists two vector functions $S$ and $S^{\ast }$ such that 
\begin{equation}
\left. 
\begin{array}{c}
S\in B\left( R\right) ,\text{ }S^{\ast }\in B^{\ast }\left( R-\delta \right)
, \\ 
\left\Vert S-S^{\ast }\right\Vert _{H}<\delta .%
\end{array}%
\right.  \label{5.7}
\end{equation}%
As to the functions $u_{0}\left( x\right) $ and $m_{0}\left( x\right) ,$ we
assume (\ref{5.02}) and also that 
\begin{equation}
\left\Vert u_{0}-u_{0}^{\ast }\right\Vert _{H^{2}\left( \Omega \right) },%
\text{ }\left\Vert m_{0}-m_{0}^{\ast }\right\Vert _{H^{2}\left( \Omega
\right) }<\delta .  \label{5.8}
\end{equation}%
Also, similarly with (\ref{3.3}) and (\ref{3.15})%
\begin{equation}
k^{\ast }\left( x\right) =\frac{2}{\left\vert \nabla u_{0}^{\ast }\left(
x\right) \right\vert ^{2}}\left[ v^{\ast }\left( x,T/2\right) +F^{\ast
}\left( x\right) \right] .  \label{5.9}
\end{equation}

For any $U\in B\left( R\right) $ and for $U^{\ast }$ denote%
\begin{equation}
\widetilde{U}=U-S,\text{ }\widetilde{U}^{\ast }=U^{\ast }-S^{\ast }.
\label{5.10}
\end{equation}%
Since $U\in B\left( R\right) $ and $U^{\ast }\in B^{\ast }\left( R-\delta
\right) ,$ then it follows from (\ref{4.5}) and (\ref{5.6}) and (\ref{5.7})
that 
\begin{equation}
\widetilde{U},\widetilde{U}^{\ast }\in \overline{B_{0}\left( 2R\right) }.
\label{5.12}
\end{equation}%
Introduce a new functional $I_{\lambda ,\beta },$ 
\begin{equation*}
I_{\lambda ,\beta }:\overline{B_{0}\left( 2R\right) }\rightarrow \mathbb{R},%
\mbox{ }I_{\lambda ,\beta }\left( V\right) =J_{\lambda ,\beta }\left(
V+S\right) ,\text{ }\forall V\in \overline{B_{0}\left( 2R\right) }
\end{equation*}%
Using triangle inequality and (\ref{5.7}), we obtain 
\begin{equation}
V+S\in \overline{B\left( 3R\right) },\mbox{ }\forall V\in \overline{%
B_{0}\left( 2R\right) }.  \label{5.13}
\end{equation}

\textbf{Theorem 5.2} (the accuracy of the minimizer and uniqueness of the
CIP). \emph{\ Suppose that conditions of Theorem 5.1 as well as conditions ( %
\ref{5.7})-( \ref{5.12}) hold. In addition, let inequality (\ref{3.1}) be
valid as well as }$\left\vert \nabla u_{0}^{\ast }\left( x\right)
\right\vert ^{2}\geq c>0$ in $\overline{\Omega }.$ \emph{Then:}

\emph{1. The functional }$I_{\lambda ,\beta }$\emph{\ has the Fr\'{e}chet
derivative }$I_{\lambda ,\beta }^{\prime }\left( V\right) \in H_{0}$\emph{\
at any point }$V\in \overline{B_{0}\left( 2R\right) }$\emph{\ and the analog
of (\ref{5.3}) holds.}

\emph{2. Let }$\overline{\lambda }_{1}=\overline{\lambda }_{1}\left( \Omega
,T,M,c,R\right) \geq 1$\emph{\ be the number of Theorem 4.3. Consider the
number }$\overline{\lambda }_{2}=\overline{\lambda }_{1}\left( \Omega
,T,M,c,3R\right) \geq \overline{\lambda }_{1}$\emph{. Then for any }$\lambda
\geq \overline{\lambda }_{2}$\emph{\ and for any choice of the
regularization parameter }$\beta \in \left[ 2\exp \left( -\lambda
T^{2}/4\right) ,1\right) $\emph{\ the functional }$I_{\lambda ,\beta }$\emph{%
\ is strongly convex on the set }$\overline{B_{0}\left( 2R\right) }$\emph{,
has unique minimizer }$V_{\min ,\lambda ,\beta }\in \overline{B_{0}\left(
2R\right) }$\emph{\ on this set, and the analog of (\ref{5.5}) holds.}

\emph{3. Let }$\sigma \in \left( 0,1\right) $\emph{\ be an arbitrary number.
Choose the number }$\gamma =\gamma \left( \rho \right) $\emph{\ as }%
\begin{equation*}
\gamma =\gamma \left( \sigma \right) =\left( \frac{\sigma }{3-\sigma }%
\right) ^{1/2}\in \left( 0,\frac{1}{\sqrt{2}}\right) .
\end{equation*}%
\emph{Choose the number }$\delta _{0}=\delta _{0}\left( \Omega
,T,M,c,R,\sigma \right) \in \left( 0,1\right) $\emph{\ so small that }%
\begin{equation*}
\ln \left\{ \delta _{0}^{-8\left( T\left( 1+\gamma \left( \sigma \right)
\right) \right) ^{-2}}\right\} \geq \emph{\ }\overline{\lambda }_{2}\geq 1.
\end{equation*}%
\emph{For any }$\delta \in \left( 0,\delta _{0}\right) $\emph{\ choose in
the functional }$I_{\lambda ,\beta }$\emph{\ parameters }$\lambda =\lambda
\left( \delta ,T,\sigma \right) \geq \overline{\lambda }_{2}$\emph{\ and }$%
\beta =\beta \left( \delta ,T,\sigma \right) $\emph{\ as}%
\begin{equation*}
\lambda \left( \delta ,T,\rho \right) =\ln \left\{ \delta ^{-8\left( T\left(
1+\gamma \left( \rho \right) \right) \right) ^{-2}}\right\} ,\text{ }\beta
=2\exp \left( -\lambda T^{2}/4\right) .
\end{equation*}

\emph{Denote }%
\begin{equation}
\overline{U}_{\min ,\lambda ,\beta }=U_{\min ,\lambda ,\beta }+S\in 
\overline{B\left( 3R\right) },  \label{5.14}
\end{equation}%
\emph{\ see (\ref{5.13}). Then the following accuracy estimates hold:}%
\begin{equation}
\hspace{-1.5cm}\left\Vert \overline{U}_{\min ,\lambda ,\beta }-U^{\ast
}\right\Vert _{\left( H^{2,1}\left( Q_{\gamma T}\right) \right) ^{4}}\leq
C_{5}\delta ^{1-\sigma },\mbox{ }\forall \delta \in \left( 0,\delta
_{0}\right) ,  \label{4.17}
\end{equation}%
\begin{equation}
\left\Vert k_{\min ,\lambda ,\beta }-k^{\ast }\right\Vert _{L_{2}\left(
\Omega \right) }\leq C_{5}\delta ^{1-\sigma },\mbox{ }\forall \delta \in
\left( 0,\delta _{0}\right) ,  \label{4.18}
\end{equation}%
\emph{where the number }$C_{5}=C_{5}\left( \Omega ,T,M,c,R,\sigma \right) >0$%
\emph{\ depends only on listed parameters. In (\ref{4.18}), the function }$%
k^{\ast }\left( x\right) $\emph{\ is computed via (\ref{5.9}), and the
function }$k_{\min ,\lambda ,\beta }\left( x\right) $\emph{\ is computed via
(\ref{3.3}) and (\ref{3.15}). }

\emph{4. Next, (\ref{5.50}), (\ref{5.7}), (\ref{5.10}), (\ref{5.14}) and (%
\ref{4.17}) imply that it is reasonable to assume that }%
\begin{equation}
\overline{U}_{\min ,\lambda ,\beta }\in \overline{B\left( R\right) }.
\label{5.15}
\end{equation}%
\emph{In the case of} \emph{(\ref{5.15}) the vector function }$\overline{U}%
_{\min ,\lambda ,\beta }$ \emph{is the unique minimizer }$\overline{U}_{\min
,\lambda ,\beta }$\ \emph{of the functional }$J_{\lambda ,\beta }$\emph{\ on
the set }$\overline{B\left( R\right) },$\emph{\ which is found in Theorem
5.1, i.e. }%
\begin{equation}
\overline{U}_{\min ,\lambda ,\beta }=U_{\min ,\lambda ,\beta },  \label{5.16}
\end{equation}%
\emph{\ and, therefore, estimate (\ref{4.17}) remains valid for }$U_{\min
,\lambda ,\beta }.$

\subsubsection{The gradient descent method}

Below in this sub-subsection all parameters are the same as the ones chosen
in Theorem 5.2. Assume now that%
\begin{equation}
U^{\ast }\in B^{\ast }\left( R/3-\delta \right) ,\mbox{ }R/3-\delta >0,
\label{4.19}
\end{equation}%
\begin{equation}
\overline{U}_{\min ,\lambda ,\beta }=U_{\min ,\lambda ,\beta }\in B\left(
R/3\right) .  \label{4.190}
\end{equation}%
As to (\ref{4.190}), also see (\ref{5.15}) and (\ref{5.16}). We construct
now the gradient descent method of the minimization of the functional $%
J_{\lambda ,\beta }.$ Consider an arbitrary vector function 
\begin{equation}
U_{0}\in B\left( R/3\right) .  \label{4.20}
\end{equation}%
Let $\xi >0$ be the step size of the gradient descent method. We construct
the following sequence of this method:%
\begin{equation}
U_{n}=U_{n-1}-\xi J_{\lambda ,\beta }^{\prime }\left( U_{n-1}\right) ,%
\mbox{
}n=1,2,...  \label{4.21}
\end{equation}%
Note that since $J_{\lambda ,\beta }^{\prime }\in H_{0}$ by Theorem 5.1,
then all vector functions $U_{n}$ have the same boundary conditions (\ref%
{3.14}), see (\ref{5.1}).

\textbf{Theorem 5.3}.\emph{\ Assume that conditions (\ref{4.19})-(\ref{4.21}%
) hold. Let all parameters of the functional }$J_{\lambda ,\beta }$\emph{\
be the same as in Theorem 5.2. Then there exists a number }$\xi _{0}\in
\left( 0,1\right) $\emph{\ such that for any }$\xi \in \left( 0,\xi
_{0}\right) $ \emph{there exists a number }$\theta =\theta \left( \xi
\right) \in \left( 0,1\right) $\emph{\ such that for all }$n\geq 1$\emph{\ } 
\begin{equation}
\hspace{-1cm}\left. 
\begin{array}{c}
U_{n}\in B\left( R\right) ,\mbox{ } \\ 
\left\Vert U_{n}-U^{\ast }\right\Vert _{H}+\left\Vert k_{n}-k^{\ast
}\right\Vert _{L_{2}\left( \Omega \right) }\leq \\ 
\leq C_{5}\delta ^{1-\sigma }+C_{5}\theta ^{n}\left\Vert U_{\min ,\lambda
,\beta }-U_{0}\right\Vert _{H},%
\end{array}%
\right.  \label{4.22}
\end{equation}%
\emph{where functions }$k_{n}\left( x\right) $\emph{\ are computed via the
direct analogs of (\ref{3.15}), and the number }$C_{5}>0$ \emph{depends on
the same parameters as ones in Theorem 5.2.}

Since $R>0$\ is an arbitrary number and in (\ref{4.20}) $U_{0}$\ is an
arbitrary point of the set $B\left( R/3\right) ,$\ then Theorem 5.3 implies
the global convergence of the gradient descent method (\ref{4.21}), see
Definition in Introduction. Formula (\ref{4.22}) provides an explicit
convergence rate.

\section{Proof of Theorem 5.1}

\label{sec:6}

Let 
\begin{equation}
U_{1}=\left( v_{1},w_{1},p_{1},q_{1}\right) \in \overline{B\left( R\right) }%
\text{ and }U_{2}=\left( v_{2},w_{2},p_{2},q_{2}\right) \in \overline{%
B\left( R\right) }  \label{6.0}
\end{equation}%
be two arbitrary points of the set $\overline{B\left( R\right) }$ defined in
(\ref{4.5}). Consider their difference 
\begin{align}
U_{2}-U_{1}=h=\left( h_{1},h_{2},h_{3},h_{4}\right) \in \overline{%
B_{0}\left( 2R\right) },  \label{6.1} \\
B_{0}\left( 2R\right) =\left\{ h\in H_{0}:\left\Vert h\right\Vert
_{H}<2R\right\} .  \label{6.2}
\end{align}%
It follows from (\ref{4.5}), (\ref{4.6}), (\ref{5.4}) and (\ref{6.0})-(\ref%
{6.2}) that we need to evaluate the differences%
\begin{equation}
\left. 
\begin{array}{c}
P_{i}=\left[ L_{i}\left(
v_{1}+h_{1},w_{1}+h_{2},p_{1}+h_{3},q_{1}+h_{4}\right) \right] ^{2}- \\ 
-\left[ L_{i}\left( v_{1},w_{1},p_{1},q_{1}\right) \right] ^{2},\text{ }%
i=1,2,3,4.%
\end{array}%
\right.  \label{6.3}
\end{equation}%
First, consider the expression $L_{1}\left(
v_{1}+h_{1},w_{1}+h_{2},p_{1}+h_{3},q_{1}+h_{4}\right) .$ By (\ref{3.10}) 
\begin{equation}
\left. 
\begin{array}{c}
L_{1}\left( v_{1}+h_{1},w_{1}+h_{2},p_{1}+h_{3},q_{1}+h_{4}\right) = \\ 
=\left( \partial _{t}+\Delta \right) \left( v_{1}+h_{1}\right) - \\ 
-2\nabla \left( v_{1}+h_{1}\right) \left( \int\limits_{T/2}^{t}\nabla \left(
v_{1}+h_{1}\right) \left( x,\tau \right) d\tau +\nabla u_{0}\left( x\right)
\right) \times \\ 
\times \left[ \left( \left( v_{1}+h_{1}\right) \left( x,t\right)
-\int\limits_{T/2}^{t}\left( w_{1}+h_{2}\right) \left( x,\tau \right) d\tau
\right) +F\left( x\right) \right] \left\vert \nabla u_{0}\left( x\right)
\right\vert ^{-2}+ \\ 
+\int\limits_{\Omega }K\left( x,y\right) \left( p_{1}+h_{3}\right) \left(
y,t\right) dy+f\left( x,t\right) \left( p_{1}+h_{3}\right) \left( x,t\right)
+ \\ 
+f_{t}\left( x,t\right) \left( \int\limits_{T/2}^{t}\left(
p_{1}+h_{3}\right) \left( x,\tau \right) d\tau +m_{0}\left( x\right) \right)
,\text{ }\left( x,t\right) \in Q_{T}.%
\end{array}%
\right.  \label{6.4}
\end{equation}%
Hence, 
\begin{equation}
\left. 
\begin{array}{c}
L_{1}\left( v_{1}+h_{1},w_{1}+h_{2},p_{1}+h_{3},q_{1}+h_{4}\right) = \\ 
=L_{1}\left( v_{1},w_{1},p_{1},q_{1}\right) +L_{1,\text{lin}}\left( h\right)
+L_{1,\text{nonlin}}\left( h\right) ,%
\end{array}%
\right.  \label{6.5}
\end{equation}%
where $L_{1,\text{lin}}$ and $L_{1,\text{nonlin}}$ are linear and nonlinear
operators respectively with respect to the vector function $h$ defined in (%
\ref{6.1}). Using (\ref{6.1}), we obtain:%
\begin{equation}
L_{1,\text{lin}}\left( h\right) =\partial _{t}h_{1}+\Delta h_{1}+\widehat{L}%
_{1,\text{lin}}\left( h,\nabla h,\int\limits_{T/2}^{t}\nabla h\left( x,\tau
\right) d\tau ,\int\limits_{T/2}^{t}h\left( x,\tau \right) d\tau ,x,t\right)
,  \label{6.6}
\end{equation}%
where $\widehat{L}_{1,\text{lin}}$ depends linearly on each of its first
four arguments. Next, using Cauchy-Schwarz inequality, Lemma 4.1, (\ref{4.04}%
), (\ref{4.4}), (\ref{4.5}), (\ref{5.02}), (\ref{6.2}) and (\ref{6.4}), we
obtain two estimates. The first estimate is:%
\begin{equation}
\left. 
\begin{array}{c}
e^{-2\lambda b^{2}}\int\limits_{Q_{T}}2\lambda ^{3/2}\left[ L_{1,\text{lin}%
}\left( h\right) +L_{1}\left( v_{1},w_{1},p_{1},q_{1}\right) \right] L_{1,%
\text{nonlin}}\left( h\right) \varphi _{\lambda }\left( x,t\right) dxdt+ \\ 
+e^{-2\lambda b^{2}}\int\limits_{Q_{T}}\lambda ^{3/2}\left( L_{1,\text{nonlin%
}}\left( h\right) \right) ^{2}\varphi _{\lambda }\left( x,t\right) dxdt\leq
K_{1}\left\Vert h\right\Vert _{H}^{2},%
\end{array}%
\right.  \label{6.7}
\end{equation}%
where the number $K_{1}=K_{1}\left( \Omega ,T,M,c,R,\lambda \right) >0$
depends only from listed parameters. The second estimate is the estimate
from the below:%
\begin{equation*}
\left. 
\begin{array}{c}
e^{-2\lambda b^{2}}\int\limits_{Q_{T}}2\lambda ^{3/2}\left[ L_{1,\text{lin}%
}\left( h\right) +L_{1}\left( v_{1},w_{1},p_{1},q_{1}\right) \right] L_{1,%
\text{nonlin}}\left( h\right) \varphi _{\lambda }\left( x,t\right) dxdt+ \\ 
+e^{-2\lambda b^{2}}\int\limits_{Q_{T}}\lambda ^{3/2}\left( L_{1,\text{nonlin%
}}\left( h\right) \right) ^{2}\varphi _{\lambda }\left( x,t\right) dxdt\geq
\\ 
\geq \left( 1/2\right) e^{-2\lambda b^{2}}\int\limits_{Q_{T}}\lambda
^{3/2}\left( \partial _{t}h_{1}+\Delta h_{1}\right) ^{2}\varphi _{\lambda
}\left( x,t\right) dxdt- \\ 
-C_{4}e^{-2\lambda b^{2}}\int\limits_{Q_{T}}\lambda ^{3/2}\left[ \left(
\nabla h_{1}\right) ^{2}+h_{1}^{2}+h_{3}^{2}+\left(
\int\limits_{T/2}^{t}\nabla h_{1}\left( x,\tau \right) d\tau \right) ^{2}%
\right] \varphi _{\lambda }\left( x,t\right) dxdt- \\ 
-C_{4}e^{-2\lambda b^{2}}\int\limits_{Q_{T}}\lambda ^{3/2}\left[ \left(
\int\limits_{T/2}^{t}h_{2}\left( x,\tau \right) d\tau \right) ^{2}+\left(
\int\limits_{T/2}^{t}h_{3}\left( x,\tau \right) d\tau \right) ^{2}\right]
\varphi _{\lambda }\left( x,t\right) dxdt%
\end{array}%
\right.
\end{equation*}%
Combining this with Lemma 4.2, we obtain the estimate, in which Volterra
integrals are not present:%
\begin{equation}
\left. 
\begin{array}{c}
e^{-2\lambda b^{2}}\int\limits_{Q_{T}}\left( L_{1,\text{lin}}\left( h\right)
\right) ^{2}\varphi _{\lambda }\left( x,t\right) dxdt+ \\ 
e^{-2\lambda b^{2}}\int\limits_{Q_{T}}2\lambda ^{3/2}\left[ L_{1,\text{lin}%
}\left( h\right) +L_{1}\left( v_{1},w_{1},p_{1},q_{1}\right) \right] L_{1,%
\text{nonlin}}\left( h\right) \varphi _{\lambda }\left( x,t\right) dxdt+ \\ 
+e^{-2\lambda b^{2}}\int\limits_{Q_{T}}\lambda ^{3/2}\left( L_{1,\text{nonlin%
}}\left( h\right) \right) ^{2}\varphi _{\lambda }\left( x,t\right) dxdt\geq
\\ 
\geq \left( 1/2\right) e^{-2\lambda b^{2}}\int\limits_{Q_{T}}\lambda
^{3/2}\left( \partial _{t}h_{1}+\Delta h_{1}\right) ^{2}\varphi _{\lambda
}\left( x,t\right) dxdt- \\ 
-C_{4}e^{-2\lambda b^{2}}\int\limits_{Q_{T}}\lambda ^{3/2}\left( \left(
\nabla h_{1}\right) ^{2}+h_{1}^{2}+h_{3}^{2}\right) \varphi _{\lambda
}\left( x,t\right) dxdt.%
\end{array}%
\right.  \label{6.8}
\end{equation}%
It is clear from (\ref{3.11})-(\ref{3.13}) that representations similar to
the one in (\ref{6.5}) are valid for $L_{i},$ $i=2,3,4,$%
\begin{equation}
\left. 
\begin{array}{c}
L_{i}\left( v_{1}+h_{1},w_{1}+h_{2},p_{1}+h_{3},q_{1}+h_{4}\right) = \\ 
=L_{i}\left( v_{1},w_{1},p_{1},q_{1}\right) +L_{i,\text{lin}}\left( h\right)
+L_{i,\text{nonlin}}\left( h\right) ,\text{ }i=2,3,4%
\end{array}%
\right.  \label{6.9}
\end{equation}%
and estimates from the above, which are similar with the one in (\ref{6.7})
are valid with constants $K_{i}.$ Denote 
\begin{equation}
L_{\text{lin}}\left( h\right) =\lambda ^{3/2}\left( L_{1,\text{lin}}\left(
h\right) +L_{2,\text{lin}}\left( h\right) \right) +L_{2,\text{lin}}\left(
h\right) +L_{4,\text{lin}}\left( h\right) .  \label{6.10}
\end{equation}%
Consider the linear functional $\widehat{J}_{\lambda ,\beta }\left( h\right)
:H_{0}\rightarrow \mathbb{R},$%
\begin{equation}
\widehat{J}_{\lambda ,\beta }\left( h\right) =e^{-2\lambda
b^{2}}\int\limits_{Q_{T}}L_{\text{lin}}\left( h\right) \varphi _{\lambda
}\left( x,t\right) dxdt+2\beta \left[ U_{1},h\right] .  \label{6.11}
\end{equation}%
By Riesz theorem there exists unique vector function $S\in H_{0}$ such that 
\begin{equation}
\widehat{J}_{\lambda ,\beta }\left( h\right) =\left[ S,h\right] ,
\label{6.12}
\end{equation}%
where $H_{0}$ is the space defined in (\ref{5.1}). Hence, it follows from (%
\ref{4.6}), (\ref{6.3}), (\ref{6.7}) and (\ref{6.9})-(\ref{6.12}) that 
\begin{equation}
\lim_{\left\Vert h\right\Vert _{H}\rightarrow 0}\frac{J_{\lambda ,\beta
}\left( U_{1}+h\right) -J_{\lambda ,\beta }\left( U_{1}\right) -\left[ S,h%
\right] }{\left\Vert h\right\Vert _{H}}=0,\text{ }\forall U_{1},U_{1}+h\in 
\overline{B\left( R\right) }.  \label{1}
\end{equation}%
Hence, 
\begin{equation}
S=J_{\lambda ,\beta }^{\prime }\left( U_{1}\right) \in H_{0}  \label{2}
\end{equation}%
is the Fr\'{e}chet derivative of the functional $J_{\lambda ,\beta }\left(
v,w,p,q\right) =J_{\lambda ,\beta }\left( U\right) $ at the point $U_{1}\in 
\overline{B\left( R\right) }.$ We omit the proof of estimate (\ref{5.3})
since it is similar with the proof of Theorem 3.1 of \cite{Bak}\emph{.}

We now come back to the estimates from the below, which are similar with the
one for $L_{1}$ in (\ref{6.8}). Consider $L_{3}$ first. Using (\ref{3.12}),
we obtain 
\begin{equation}
\left. 
\begin{array}{c}
e^{-2\lambda b^{2}}\int\limits_{Q_{T}}\lambda ^{3/2}\left( L_{3,\text{lin}%
}\left( h\right) \right) ^{2}\varphi _{\lambda }\left( x,t\right) dxdt+ \\ 
+e^{-2\lambda b^{2}}\times \\ 
\times \int\limits_{Q_{T}}2\lambda ^{3/2}\left[ L_{3,\text{lin}}\left(
h\right) +L_{3}\left( v_{1},w_{1},p_{1},q_{1}\right) \right] L_{3,\text{%
nonlin}}\left( h\right) \varphi _{\lambda }\left( x,t\right) dxdt+ \\ 
+e^{-2\lambda b^{2}}\int\limits_{Q_{T}}\lambda ^{3/2}\left( L_{3,\text{nonlin%
}}\left( h\right) \right) ^{2}\varphi _{\lambda }\left( x,t\right) dxdt\geq
\\ 
\geq \left( 1/2\right) e^{-2\lambda b^{2}}\int\limits_{Q_{T}}\lambda
^{3/2}\left( \partial _{t}h_{2}+\Delta h_{2}\right) ^{2}\varphi _{\lambda
}\left( x,t\right) dxdt- \\ 
-C_{4}e^{-2\lambda b^{2}}\int\limits_{Q_{T}}\lambda ^{3/2}\left( \left(
\nabla h_{1}\right) ^{2}+\left( \nabla h_{2}\right)
^{2}+h_{1}^{2}+h_{2}^{2}+h_{3}^{2}+h_{4}^{2}\right) \varphi _{\lambda
}\left( x,t\right) dxdt.%
\end{array}%
\right.  \label{6.13}
\end{equation}%
Similarly consider $L_{2},$ using (\ref{3.11}),%
\begin{equation}
\left. 
\begin{array}{c}
e^{-2\lambda b^{2}}\int\limits_{Q_{T}}\left( L_{2,\text{lin}}\left( h\right)
\right) ^{2}\varphi _{\lambda }\left( x,t\right) dxdt+ \\ 
+e^{-2\lambda b^{2}}\int\limits_{Q_{T}}2\left[ L_{2,\text{lin}}\left(
h\right) +L_{2}\left( v_{1},w_{1},p_{1},q_{1}\right) \right] L_{2,\text{%
nonlin}}\left( h\right) \varphi _{\lambda }\left( x,t\right) dxdt+ \\ 
+e^{-2\lambda b^{2}}\int\limits_{Q_{T}}\left( L_{2,\text{nonlin}}\left(
h\right) \right) ^{2}\varphi _{\lambda }\left( x,t\right) dxdt\geq \\ 
\geq \left( 1/2\right) e^{-2\lambda b^{2}}\int\limits_{Q_{T}}\left( \partial
_{t}h_{3}-\Delta h_{3}\right) ^{2}\varphi _{\lambda }\left( x,t\right) dxdt-
\\ 
-C_{4}e^{-2\lambda b^{2}}\int\limits_{Q_{T}}\left( \left( \nabla
h_{1}\right) ^{2}+\left( \nabla h_{2}\right) ^{2}+\left( \nabla h_{3}\right)
^{2}+\left( \nabla h_{4}\right) ^{2}\right) \varphi _{\lambda }\left(
x,t\right) dxdt- \\ 
-C_{4}e^{-2\lambda b^{2}}\int\limits_{Q_{T}}\left(
h_{1}^{2}+h_{2}^{2}+h_{3}^{2}+h_{4}^{2}\right) \varphi _{\lambda }\left(
x,t\right) dxdt- \\ 
-C_{4}e^{-2\lambda b^{2}}\int\limits_{Q_{T}}\left( \Delta h_{1}\right)
^{2}\varphi _{\lambda }\left( x,t\right) dxdt.%
\end{array}%
\right.  \label{6.14}
\end{equation}%
Similarly, consider $L_{4}.$ Using (\ref{3.13}), we obtain%
\begin{equation}
\left. 
\begin{array}{c}
e^{-2\lambda b^{2}}\int\limits_{Q_{T}}\left( L_{4,\text{lin}}\left( h\right)
\right) ^{2}\varphi _{\lambda }\left( x,t\right) dxdt+ \\ 
+e^{-2\lambda b^{2}}\int\limits_{Q_{T}}2\left[ L_{4,\text{lin}}\left(
h\right) +L_{4}\left( v_{1},w_{1},p_{1},q_{1}\right) \right] L_{4,\text{%
nonlin}}\left( h\right) \varphi _{\lambda }\left( x,t\right) dxdt+ \\ 
+e^{-2\lambda b^{2}}\int\limits_{Q_{T}}\left( L_{4,\text{nonlin}}\left(
h\right) \right) ^{2}\varphi _{\lambda }\left( x,t\right) dxdt\geq \\ 
\geq \left( 1/2\right) e^{-2\lambda b^{2}}\int\limits_{Q_{T}}\left( \partial
_{t}h_{4}-\Delta h_{4}\right) ^{2}\varphi _{\lambda }\left( x,t\right) dxdt-
\\ 
-C_{4}e^{-2\lambda b^{2}}\int\limits_{Q_{T}}\left( \left( \nabla
h_{1}\right) ^{2}+\left( \nabla h_{2}\right) ^{2}+\left( \nabla h_{3}\right)
^{2}+\left( \nabla h_{4}\right) ^{2}\right) \varphi _{\lambda }\left(
x,t\right) dxdt- \\ 
-C_{4}e^{-2\lambda b^{2}}\int\limits_{Q_{T}}\left(
h_{1}^{2}+h_{2}^{2}+h_{3}^{2}+h_{4}^{2}\right) \varphi _{\lambda }\left(
x,t\right) dxdt- \\ 
-C_{4}e^{-2\lambda b^{2}}\int\limits_{Q_{T}}\left( \Delta h_{1}\right)
^{2}\varphi _{\lambda }\left( x,t\right) dxdt-C_{4}e^{-2\lambda
b^{2}}\int\limits_{Q_{T}}\left( \Delta h_{2}\right) ^{2}\varphi _{\lambda
}\left( x,t\right) dxdt.%
\end{array}%
\right.  \label{6.15}
\end{equation}

A significant difficulty for the further analysis is due to the presence of
the terms with $\left( \Delta h_{1}\right) ^{2}$ and $\left( \Delta
h_{2}\right) ^{2}$ in (\ref{6.14}) and (\ref{6.15}). This is exactly a
reflection of the second difficulty of working with the MFGS, which was
mentioned in subsection 2.2. It is because of the presence of these two
terms why the multiplier $\lambda ^{3/2}$ in (\ref{6.8}) and (\ref{6.13}).

Recall that we need to prove the strong convexity estimate (\ref{5.4}).
Thus, using (\ref{1})-(\ref{6.15}), we obtain%
\begin{equation}
\left. 
\begin{array}{c}
J_{\lambda ,\beta }\left( U_{1}+h\right) -J_{\lambda ,\beta }\left(
U_{1}\right) -\left[ J_{\lambda ,\beta }^{\prime }\left( U_{1}\right) ,h%
\right] \geq \\ 
\geq \left( 1/2\right) e^{-2\lambda b^{2}}\int\limits_{Q_{T}}\lambda
^{3/2}\left( \partial _{t}h_{1}+\Delta h_{1}\right) ^{2}\varphi _{\lambda
}\left( x,t\right) dxdt+ \\ 
+\left( 1/2\right) e^{-2\lambda b^{2}}\int\limits_{Q_{T}}\lambda
^{3/2}\left( \partial _{t}h_{2}+\Delta h_{2}\right) ^{2}\varphi _{\lambda
}\left( x,t\right) dxdt+ \\ 
+\left( 1/2\right) e^{-2\lambda b^{2}}\int\limits_{Q_{T}}\left( \partial
_{t}h_{3}-\Delta h_{3}\right) ^{2}\varphi _{\lambda }\left( x,t\right) dxdt+
\\ 
+\left( 1/2\right) e^{-2\lambda b^{2}}\int\limits_{Q_{T}}\left( \partial
_{t}h_{4}-\Delta h_{4}\right) ^{2}\varphi _{\lambda }\left( x,t\right) dxdt-
\\ 
-C_{4}e^{-2\lambda b^{2}}\int\limits_{Q_{T}}\lambda ^{3/2}\left( \left(
\nabla h_{1}\right) ^{2}+\left( \nabla h_{2}\right) ^{2}\right) \varphi
_{\lambda }\left( x,t\right) dxdt- \\ 
-C_{4}e^{-2\lambda b^{2}}\int\limits_{Q_{T}}\lambda ^{3/2}\left(
h_{1}^{2}+h_{2}^{2}+h_{3}^{2}+h_{4}^{2}\right) \varphi _{\lambda }\left(
x,t\right) dxdt- \\ 
-C_{4}e^{-2\lambda b^{2}}\int\limits_{Q_{T}}\left( \left( \nabla
h_{3}\right) ^{2}+\left( \nabla h_{4}\right) ^{2}\right) \varphi _{\lambda
}\left( x,t\right) dxdt- \\ 
-C_{4}e^{-2\lambda b^{2}}\int\limits_{Q_{T}}\left( \left( \Delta
h_{1}\right) ^{2}+\left( \Delta h_{2}\right) ^{2}\right) \varphi _{\lambda
}\left( x,t\right) dxdt,\text{ }\forall \lambda \geq \overline{\lambda }_{1},%
\end{array}%
\right.  \label{6.16}
\end{equation}%
where $\overline{\lambda }_{1}\geq \overline{\lambda }_{0}$ is a
sufficiently large number depending only on parameters listed in the
formulation of this theorem, and $\overline{\lambda }_{0}\geq 1$ is the
number of Lemma 4.3. We will specify more $\overline{\lambda }_{1}$ later.

Apply now the Carleman estimate of Lemma 4.3 to lines 2 and 3 of (\ref{6.16}%
). To estimate terms with $\left\Vert h_{k}\left( x,T\right) \right\Vert
_{H^{1}\left( \Omega \right) }^{2},\left\Vert h_{k}\left( x,0\right)
\right\Vert _{H^{1}\left( \Omega \right) }^{2},$ we use trace theorem. In
addition, we use the second line of (\ref{4.2}). We obtain%
\begin{equation}
\left. 
\begin{array}{c}
\left( 1/2\right) e^{-2\lambda b^{2}}\int\limits_{Q_{T}}\lambda ^{3/2}\left(
\partial _{t}h_{1}+\Delta h_{1}\right) ^{2}\varphi _{\lambda }\left(
x,t\right) dxdt+ \\ 
+\left( 1/2\right) e^{-2\lambda b^{2}}\int\limits_{Q_{T}}\lambda
^{3/2}\left( \partial _{t}h_{2}+\Delta h_{2}\right) ^{2}\varphi _{\lambda
}\left( x,t\right) dxdt\geq \\ 
\geq C_{3}e^{-2\lambda b^{2}}\sum\limits_{k=1}^{2}\sqrt{\lambda }%
\int\limits_{Q_{T}}\left( \left( \partial _{t}h_{k}\right)
^{2}+\sum\limits_{i,j=1}^{n}\left( \partial _{x_{i}}\partial
_{x_{j}}h_{k}\right) ^{2}\right) \varphi _{\lambda }\left( x,t\right) dxdt+
\\ 
+C_{3}e^{-2\lambda b^{2}}\sum\limits_{k=1}^{2}\lambda
^{5/2}\int\limits_{Q_{T}}\left( \left( \nabla h_{k}\right) ^{2}+\lambda
^{2}h_{k}^{2}\right) \varphi _{\lambda }\left( x,t\right) dxdt- \\ 
-C_{4}\lambda ^{3/2}e^{-\lambda T^{2}/2}\sum\limits_{k=1}^{2}\left\Vert
h_{k}\right\Vert _{H^{k_{n}}\left( Q_{T}\right) }^{2},\text{ }\forall
\lambda \geq \overline{\lambda }_{1}.%
\end{array}%
\right.  \label{6.17}
\end{equation}%
Choose $\overline{\lambda }_{1}\geq 1$ so large that $C_{3}\sqrt{\overline{%
\lambda }_{1}}>2C_{4}.$ Then sum up (\ref{6.17}) with lines number 6-9 of (%
\ref{6.16}). We obtain%
\begin{equation}
\left. 
\begin{array}{c}
J_{\lambda ,\beta }\left( U_{1}+h\right) -J_{\lambda ,\beta }\left(
U_{1}\right) -\left[ J_{\lambda ,\beta }^{\prime }\left( U_{1}\right) ,h%
\right] \geq \\ 
\geq C_{4}e^{-2\lambda b^{2}}\sum\limits_{k=1}^{2}\sqrt{\lambda }%
\int\limits_{Q_{T}}\left( \left( \partial _{t}h_{k}\right)
^{2}+\sum\limits_{i,j=1}^{n}\left( \partial _{x_{i}}\partial
_{x_{j}}h_{k}\right) ^{2}\right) \varphi _{\lambda }\left( x,t\right) dxdt+
\\ 
+C_{4}e^{-2\lambda b^{2}}\sum\limits_{k=1}^{2}\lambda
^{5/2}\int\limits_{Q_{T}}\left( \left( \nabla h_{k}\right) ^{2}+\lambda
^{2}h_{k}^{2}\right) \varphi _{\lambda }\left( x,t\right) dxdt+ \\ 
+\left( 1/2\right) e^{-2\lambda b^{2}}\int\limits_{Q_{T}}\left( \partial
_{t}h_{3}-\Delta h_{3}\right) ^{2}\varphi _{\lambda }\left( x,t\right) dxdt+
\\ 
+\left( 1/2\right) e^{-2\lambda b^{2}}\int\limits_{Q_{T}}\left( \partial
_{t}h_{4}-\Delta h_{4}\right) ^{2}\varphi _{\lambda }\left( x,t\right) dxdt-
\\ 
-C_{4}e^{-2\lambda b^{2}}\int\limits_{Q_{T}}\left( \left( \nabla
h_{3}\right) ^{2}+\left( \nabla h_{4}\right) ^{2}+\lambda
^{3/2}h_{3}^{2}+\lambda ^{3/2}h_{4}^{2}\right) \varphi _{\lambda }\left(
x,t\right) dxdt- \\ 
-C_{4}\lambda ^{3/2}e^{-\lambda T^{2}/2}\sum\limits_{k=1}^{2}\left\Vert
h_{k}\right\Vert _{H^{k_{n}}\left( Q_{T}\right) }^{2}+\beta \left\Vert
h\right\Vert _{\left( H^{k_{n}}\left( Q_{T}\right) \right) ^{4}}^{2},\text{ }%
\forall \lambda \geq \overline{\lambda }_{1}.%
\end{array}%
\right.  \label{6.18}
\end{equation}%
To absorb the negative integral in the sixth line of (\ref{6.18}), we apply
the Carleman estimate of Lemma 4.3 to the fourth and fifth lines of (\ref%
{6.18}). Increasing, if necessary $\overline{\lambda }_{1}$ and taking into
account (\ref{4.2}), we obtain%
\begin{equation}
\left. 
\begin{array}{c}
J_{\lambda ,\beta }\left( U_{1}+h\right) -J_{\lambda ,\beta }\left(
U_{1}\right) -\left[ J_{\lambda ,\beta }^{\prime }\left( U_{1}\right) ,h%
\right] \geq \\ 
\geq C_{4}\exp \left[ -\lambda \left( \gamma T\right) ^{2}/2\right]
\left\Vert h\right\Vert _{\left( H^{2,1}\left( Q_{\gamma ,T}\right) \right)
^{4}}^{2}+ \\ 
+\left( \beta -C_{4}\lambda ^{3/2}\exp \left( -\lambda T^{2}/2\right)
\right) \left\Vert h\right\Vert _{H}^{2},\text{ }\forall \lambda \geq 
\overline{\lambda }_{1},\text{ }\forall \gamma \in \left( 0,1\right] .%
\end{array}%
\right.  \label{6.19}
\end{equation}%
Since the regularization parameter $\beta \in \left[ 2\exp \left( -\lambda
T^{2}/4\right) ,1\right) ,$ then (\ref{6.19}) implies the target estimate (%
\ref{5.4}) of this theorem. As soon as (\ref{5.4}) is established, existence
and uniqueness of the minimizer $U_{\min ,\lambda ,\beta }$ of the
functional $J_{\lambda ,\beta }\left( U\right) $ on the set $\overline{%
B\left( R\right) }$ as well as inequality (\ref{5.4}) follow from a
combination of Lemma 2.1 with Theorem 2.1 of \cite{Bak}\emph{. }$\square $

\section{Numerical Studies}

\label{sec:7}

In this section we describe our numerical studies of the Minimization
Problem formulated in section 4.2. First, we need to figure out how to
generate the input data (\ref{2.3}) for our CIP, and then we need to carry
out numerical tests for this inverse problem. Our data generation procedure
is the same as the one in \cite{MFG8}. For the convenience of the reader,
this procedure is described in subsection \ref{sec:8.1}. Numerical tests are
described in subsection \ref{sec:8.2}.

\subsection{Numerical generation of input data (\protect\ref{2.3})}

\label{sec:8.1}

First, we choose the coefficient $\widehat{k}\left( x\right) $, which we
want to reconstruct. Next, we choose a sufficiently smooth function $%
\widehat{u}(x,t)$ such that condition (\ref{3.1}) is satisfied. Next we
solve the following initial boundary value problem for the FP equation in (%
\ref{2.2}): 
\begin{equation}
\left. 
\begin{array}{c}
\widehat{m}_{t}(x,t)-\Delta \widehat{m}(x,t)-\text{div}(\widehat{k}(x)%
\widehat{m}(x,t)\nabla \widehat{u}(x,t))=0,\quad \left( x,t\right) \in Q_{T},
\\ 
\widehat{m}(x,0)=1,\quad x\in \Omega , \\ 
\widehat{m}\mid _{S_{T}}=\widehat{h}_{m}\left( x,t\right) ,\text{ }\left(
x,t\right) \in S_{T},%
\end{array}%
\right.  \label{8.101}
\end{equation}%
where $\widehat{h}_{m}\left( x,t\right) $ is the function of our choice with 
$\widehat{h}_{m}\left( x,t\right) =1$. We solve the initial boundary value
problem (\ref{8.101}) via the Finite Difference Method.

Given functions $\widehat{k}\left( x\right) ,\widehat{u}(x,t),\widehat{m}%
(x,t),$ we now need to ensure that the function $\widehat{u}(x,t)$ satisfies
the HJB equation in MFGS (\ref{2.2}). To do this, we define the function $%
f(x,t)$ in (\ref{2.2}) in a special way. Assume that%
\begin{equation}
\widehat{m}(x,t)\neq 0\text{ in }\overline{Q}_{T}.  \label{7.1}
\end{equation}
It follows from the maximum principle for parabolic equations that (\ref{7.1}%
) can always be achieved by a proper choice of the function $\widehat{h}%
_{m}\left( x,t\right) $ in (\ref{8.101}). Using the HJB equation in (\ref%
{2.2}), we define the function $f\left( x,t\right) $ as: 
\begin{equation}
\begin{split}
& f(x,t)=-\frac{\widehat{u}_{t}(x,t)+\Delta \widehat{u}(x,t){-\widehat{k}(}x{%
)(\nabla \widehat{u}(x,t))^{2}/2}+\int\limits_{\Omega }K\left( x,y\right) 
\widehat{m}(y,t)dy}{\widehat{m}(x,t)}, \\
& \hspace{2cm}\left( x,t\right) \in Q_{T}.
\end{split}
\label{8.104}
\end{equation}%
Hence, it follows from (\ref{8.101}) and (\ref{8.104}) that the so chosen
pair of functions $\left( \widehat{u},\widehat{m}\right) (x,t)$ satisfies
MFGS (\ref{2.2}).

This is our numerical method for the generation of a solution $\left( 
\widehat{u},\widehat{m}\right) (x,t)$ of MFGS (\ref{2.2}). Next, the values
of functions $\widehat{u}\left( x,t\right) $ and $\widehat{m}\left(
x,t\right) $ at $t=T/2$, on $S_{T}$ as well as the values of derivatives $%
\widehat{u}_{x_{1}}\left( x,t\right) $ and $\widehat{m}_{x_{1}}\left(
x,t\right) $ on $\Psi _{1T}^{+}$ generate the target data (\ref{2.3}) for
our CIP.

\subsection{Numerical arrangements}

\label{sec:8.2}

Below $x=\left( x_{1},x_{2}\right) .$ We have conducted numerical studies in
the 2D case with the domain: 
\begin{equation}
\Omega =\left( 1,2\right) \times \left( 1,2\right) ,\quad T=1.  \label{8.201}
\end{equation}%
To generate the input data (\ref{2.3}), we set:

\begin{equation}
\widehat{u}(x,t)=x_{1}^{2}x_{2}^{2}(1+t),\quad \widehat{h}_{m}\left(
x,t\right) =1+x_{1}x_{2}t.  \label{7.2}
\end{equation}

As to the kernel $K\left( x,y\right) $ in the global interaction term with $%
x=(x_{1},x_{2}),y=(y_{1},y_{2})$, we have taken it as in (\ref{2.4}), 
\begin{equation}
K\left( x,y\right) =K\left( x_{1},x_{2},y_{1},y_{2}\right) =\delta \left(
x_{1}-y_{1}\right) \exp \left( -\frac{\left( x_{2}-y_{2}\right) ^{2}}{\sigma
^{2}}\right) ,\quad \sigma =0.2.  \label{8.202}
\end{equation}%
It was noted in the beginning of subsection 2.3 that the form (\ref{8.202})
is close to the one proposed in \cite[formula (2.7)]{Ding}.

In all numerical tests below the target coefficient $k(x)$ in (\ref{2.2}),
to be reconstructed, is taken as 
\begin{equation}
k(x)=\left\{ 
\begin{array}{cc}
c_{a}=const.>1, & \text{inside the tested inclusion,} \\ 
1, & \text{outside the tested inclusion.}%
\end{array}%
\right.  \label{8.203}
\end{equation}%
Naturally, to avoid singularities of the solutions of forward problems (\ref%
{8.101}) in the above data generation procedure, we slightly smooth out $%
k(x) $ near boundaries of our tested inclusions. So that the resulting
function $k(x)\in C^{1}\left( \overline{\Omega }\right) .$Then we set: 
\begin{align}
\mbox{correct inclusion/background contrast}&=\frac{c_{a}}{1},  \label{8.204}
\\
\text{computed inclusion/background contrast}&=\frac{\max_{\text{inclusion}%
}\left( k_{\text{comp}}(x)\right) }{1}.  \label{8.2040}
\end{align}%
In the numerical tests below, we take $c_{a}=2,4,8$, and the inclusions with
the shapes of the letters `$A$', `$\Omega $' and `$SZ$'. As soon as choices (%
\ref{8.201})-(\ref{8.203}) are in place, we proceed with the numerical
generation of the input data (\ref{2.3}), as described in subsection \ref%
{sec:8.1}.

\textbf{Remark 7.1.} \emph{We choose shapes of letters for our tested
inclusions to demonstrate robustness of our numerical method. Indeed, these
are complex non-convex shapes with voids. On the other hand, since our CIP
is a quite challenging one in its own right, we are not concerned with
accurate reconstructions of edges of inclusions: it is sufficient for us, as
long as reconstructed shapes are \textquotedblleft generally" rather
accurate ones, and also as long as correct contrasts in (\ref{8.204}) are
close to the computed ones in (\ref{8.2040}).}

The regularization parameter $\beta $ in functional (\ref{4.6}) and the
parameter $\lambda $ in the Carleman Weight Function $\varphi _{\lambda }$
of this functional were: 
\begin{equation}
\beta =0.001,\quad \lambda =3,  \label{8.205}
\end{equation}%
see (\ref{4.1}) for the formula for the function $\varphi _{\lambda }.$ The
choice of $\lambda $ is a delicate one, and it is described in Test 1 below.

To solve the forward problem for data generation, we have used the spatial
mesh sizes $1/160\times 1/160$ and the temporal mesh step size $1/320$. In
the computations of the Minimization Problem, the spatial mesh sizes were $%
1/20\times 1/20$ and the temporal mesh step size was $1/10$. We have solved
problem (\ref{8.101}) by the classic implicit scheme. To solve the
Minimization Problem, we have written both differential operators and the
norm $\left\Vert \cdot \right\Vert _{H^{k_{n}}\left( Q_{T}\right) }$ in (\ref%
{4.6}) in the discrete forms of finite differences and then minimized the
functional $J_{\lambda ,\beta }\left( v,w,p,q\right) $ with respect to the
values of functions $v,w,p,q$ at those grid points. As soon as its minimizer 
$\left( v_{\min },w_{\min },p_{\min },q_{\min }\right) $ is found, the
computed target coefficient $k_{\mbox{comp}}\left( x\right) $ is found via
an obvious analog of (\ref{3.15}). Now, even though by (\ref{4.04}) we
should have norms $\left\Vert \cdot \right\Vert _{H^{k_{n}}\left(
Q_{T}\right) }=\left\Vert \cdot \right\Vert _{H^{3}\left( Q_{T}\right) }$ in
(\ref{4.6}), we have still used $H^{2}\left( Q_{T}\right) -$norms in (\ref%
{4.6}) since $H^{3}\left( Q_{T}\right) -$norms are harder to work on
computationally. We conjecture that $H^{2}\left( Q_{T}\right) -$norms work
well computationally in our tests since all norms in any finite dimensional
space are equivalent, and we have relatively small number of grid points,
which effectively means a relatively low dimensions of our spaces of
discrete functions.

To guarantee that the solution of the problem of the minimization of the
functional $J_{\lambda ,\beta }\left( v,w,p,q\right) $ in (\ref{4.6})
satisfies the boundary conditions (\ref{3.14}), we adopt the Matlab's
built-in optimization toolbox \textbf{fmincon} to minimize the discretized
form of this functional. The iterations of \textbf{fmincon} stop when the
condition 
\begin{equation}
|\nabla J_{\lambda ,\beta }\left( v,w,p,q\right) |<10^{-2}  \label{7.3}
\end{equation}%
is met. This stopping criterion is justified in Test 1 below.

It is easy to deal with the Dirichlet boundary conditions in (\ref{3.14})
for $v,w,p,q$ by just using the equations in the first column of (\ref{3.14}%
) on the discrete grid points on boundary $S_{T}$ in all iterations of 
\textbf{fmincon}. To exhibit the process for dealing with the Neumann
boundary conditions in (\ref{3.14}) in the iterations of \textbf{fmincon},
we denote the discrete points along the $x_{1}$-direction as 
\begin{equation}
x_{1,i}=1+ih_{x_{1}},\quad i=0,1,\cdots ,N_{x_{1}},\quad
N_{x_{1}}=1/h_{x},\quad h_{x_{1}}=1/20.  \label{8.206}
\end{equation}%
The finite difference method is adopted to numerically approximate the $%
x_{1} $-derivations in the second column of (\ref{3.14}) on the part $\Psi
_{1T}^{+}$ of the lateral boundary $S_{T}$. Then, keeping in mind those
Neumann boundary conditions (\ref{3.14}), the discrete functions $v,w,p,q$
in all iterations of \textbf{fmincon} should satisfy 
\begin{equation}
\left. 
\begin{array}{c}
-4v(x_{N_{x_{1}}-1},x_{2},t)+v(x_{N_{x_{1}}-2},x_{2},t)=2h_{x_{1}}\partial
_{t}g_{1}(2,x_{2},t)-3\partial _{t}g_{0}(2,x_{2},t), \\ 
-4w(x_{N_{x}-1},x_{2},t)+w(x_{N_{x}-2},x_{2},t)=2h_{x_{1}}\partial
_{t}^{2}g_{1}(2,x_{2},t)-3\partial _{t}^{2}g_{0}(2,x_{2},t), \\ 
-4p(x_{N_{x}-1},x_{2},t)+p(x_{N_{x}-2},x_{2},t)=2h_{x_{1}}\partial
_{t}p_{1}(2,x_{2},t)-3\partial _{t}p_{0}(2,x_{2},t), \\ 
-4q(x_{N_{x}-1},x_{2},t)+q(x_{N_{x}-2},x_{2},t)=2h_{x_{1}}\partial
_{t}^{2}p_{1}(2,x_{2},t)-3\partial _{t}^{2}p_{0}(2,x_{2},t).%
\end{array}%
\right. \hspace{-0.2cm}  \label{8.207}
\end{equation}%
We note that the right part of formula \eqref{8.207} also contains the
Dirichlet boundary conditions in (\ref{3.14}) on the boundary $\Psi
_{1T}^{+}\subset S_{T}$ with $x_{1}=2$ as 
\begin{equation}
\begin{split}
v(x_{N_{x_{1}}},x_{2},t)& \mid _{\Psi _{1T}^{+}}=\partial
_{t}g_{0}(2,x_{2},t),\quad w(x_{N_{x_{1}}},x_{2},t)\mid _{\Psi
_{1T}^{+}}=\partial _{t}^{2}g_{0}(2,x_{2},t), \\
\quad p(x_{N_{x_{1}}},x_{2},t)& \mid _{\Psi _{1T}^{+}}=\partial
_{t}p_{0}(2,x_{2},t),\quad q(x_{N_{x_{1}}},x_{2},t)\mid _{\Psi
_{1T}^{+}}=\partial _{t}^{2}p_{0}(2,x_{2},t).
\end{split}
\label{8.208}
\end{equation}

The starting point $\left( v^{\left( 0\right) },w^{\left( 0\right)
},p^{\left( 0\right) },q^{\left( 0\right) }\right) \left( x,t\right) $ of
iterations of \textbf{fmincon} was chosen as: 
\begin{equation}
v^{(0)}(x,t)=w^{(0)}(x,t)=p^{(0)}(x,t)=q^{(0)}(x,t)=0.  \label{8.209}
\end{equation}%
Although it follows from (\ref{8.209}) that the starting point $\left(
v^{\left( 0\right) },w^{\left( 0\right) },p^{\left( 0\right) },q^{\left(
0\right) }\right) $ does not satisfy the boundary conditions in (\ref{3.14}%
), still (\ref{8.207}) implies that boundary conditions \eqref{3.14} are
satisfied on all other iterations of \textbf{fmincon}. In the procedure of 
\textbf{fmincon}, the Dirichlet boundary conditions are ensured by making
the values of functions $v,w,p,q$ on the discrete grid points on boundary $%
S_{T}$ to satisfy the first column of (\ref{3.14}), and the Neumann boundary
conditions are ensured by (\ref{8.207}). Then the the discrete functions $%
v,w,p,q$ on the discrete grid points satisfy both Dirichlet and Neumann
boundary conditions (\ref{3.14}).

We introduce the random noise in the time dependent boundary input data in (%
\ref{3.2}) as follows: 
\begin{equation}
\begin{split}
g_{0}^{\xi }\left( x,t\right) =g_{0}\left( x,t\right) \left( 1+\delta \xi
_{g_{0}}\left( t\right) \right) ,& \quad p_{0}^{\xi }\left( x,t\right)
=p_{0}\left( x,t\right) \left( 1+\delta \xi _{p_{0}}\left( t\right) \right) ,
\\
g_{1}^{\xi }\left( x,t\right) =g_{1}\left( x,t\right) \left( 1+\delta \xi
_{g_{1}}\left( t\right) \right) ,& \quad p_{1}^{\xi }\left( x,t\right)
=p_{1}\left( x,t\right) \left( 1+\delta \xi _{p_{1}}\left( t\right) \right) ,
\end{split}
\label{8.210}
\end{equation}%
where $\xi _{g_{0}}, \xi _{p_{0}}, \xi _{g_{1}},\xi _{p_{1}}$ are the
uniformly distributed random variables in the interval $[-1,1]$ depending on
the point $t\in \lbrack 0,T]$. In (\ref{8.210}) $\delta =0.03,0.05$, which
correspond to the $3\%$ and $5\%$ noise levels respectively. The
reconstruction from the noisy data is denoted as $k^{\xi }(x)$ given by the
following analog of (\ref{3.2}): 
\begin{equation}
\left. 
\begin{array}{c}
k^{\xi }\left( x\right) =\frac{2}{\left\vert \nabla u_{0}\left( x\right)
\right\vert ^{2}}\left[ v_{\min }^{\xi }\left( x,T/2\right) +F\left(
x\right) \right] , \\ 
F\left( x\right) =\Delta u_{0}\left( x\right) +\int\nolimits_{\Omega
}K\left( x,y\right) m_{0}\left( y\right) dy+f\left( x,T/2\right) m_{0}\left(
x\right) ,%
\end{array}%
\right.  \label{8.211}
\end{equation}%
where the subscript $\xi $ means that these functions correspond to the
noisy data. Since we deal with first and second $t-$derivatives of noisy
functions $g_{0}^{\xi }\left( x,t\right) $, $p_{0}^{\xi }\left( x,t\right) $%
, $g_{1}^{\xi }\left( x,t\right) $, $p_{1}^{\xi }\left( x,t\right) $, we
have to design a numerical method to differentiate the noisy data. First, we
use the natural cubic splines to approximate the noisy input data (\ref%
{8.210}). Next, we use the derivatives of those splines to approximate the
derivatives of corresponding noisy observation data. We generate the
corresponding cubic splines in the temporal space with the temporal mesh
grid size of $1/10$, and then calculate the derivatives to approximate the
first and second derivatives with respect to $t$.

\subsection{Numerical tests}

\label{sec:8.3}

\begin{figure}[tbph]
\centering
\includegraphics[width = 4.5in]{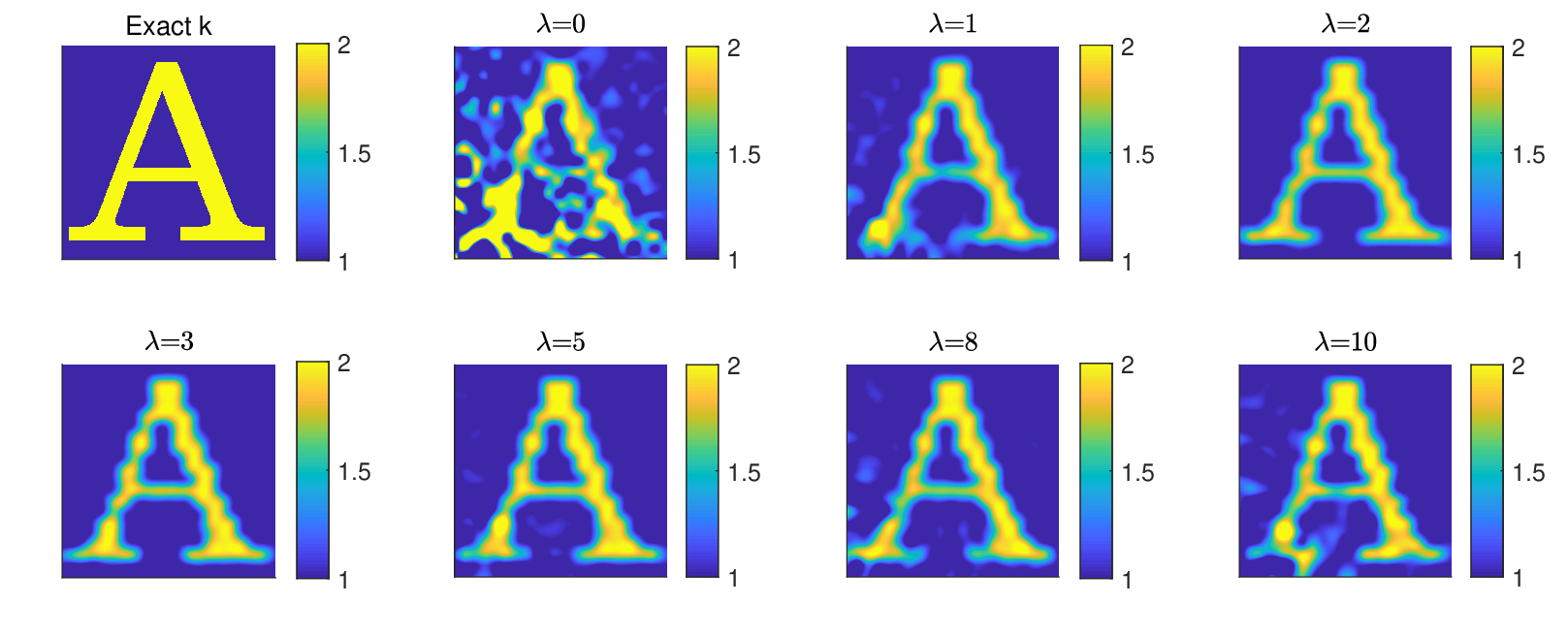}
\caption{Test 1. The reconstructed coefficient $k(x)$, where the function $%
k(x)$ is given in \eqref{8.203} with $c_{a}=2$ inside of the letter `A'. We
test different values of $\protect\lambda $.}
\label{plot_re_A2_diff_lambda}
\end{figure}

\textbf{Test 1.} We test the case when the inclusion in (\ref{8.203}) has
the shape of the letter `$A$' with $c_{a}=2$. We use this test as a
reference one to figure out an optimal value of the parameter $\lambda $.
The result is displayed in Figure \ref{plot_re_A2_diff_lambda}. We observe
that the images have a low quality for $\lambda =0,1$. Then the quality is
improved with $\lambda =2,3,5$, and the reconstruction quality deteriorates
at $\lambda =10$. On the other hand, the image is accurate at $\lambda =3,$
including the accurate reconstruction of the inclusion/background contrast (%
\ref{8.2040}).

\textbf{Conclusion:} Hence, we choose $\lambda =3$ as the optimal value and
we use this value in other tests, see (\ref{8.205}).

\begin{figure}[tbph]
\centering
\includegraphics[width = 4in]{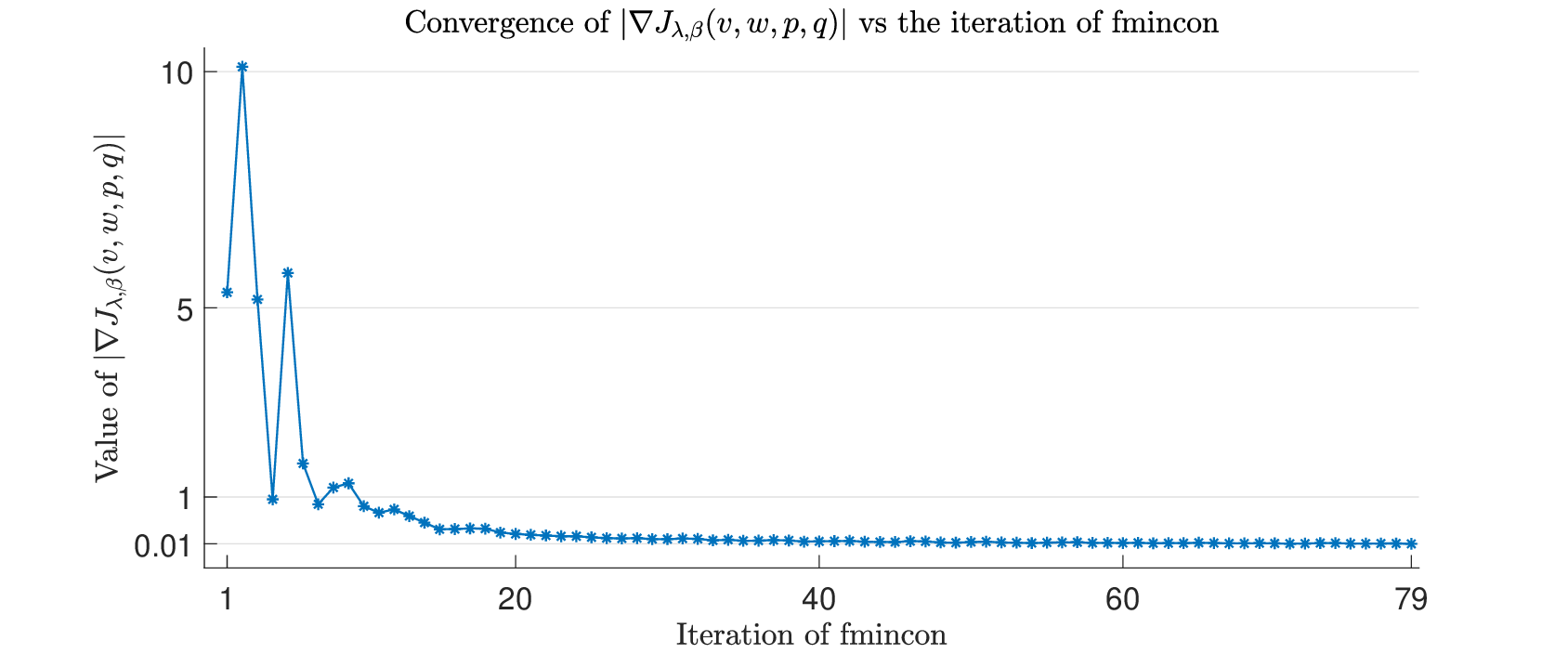}
\caption{Test 1. The convergence of $|\protect\nabla J_{\protect\lambda ,%
\protect\beta }\left( v, w, p, q \right) |$ in the iteration of \textbf{%
fmincon} with $\protect\lambda=3$ and the function $k(x)$ is given in 
\eqref{8.203} with $c_{a}=2$ inside of the letter `A'.}
\label{plot_grad}
\end{figure}

With $\lambda =3$, we display now in Figure \ref{plot_grad} the convergence
behavior of $|\nabla J_{\lambda ,\beta }$ $( v$,$w$,$p$,$q) |$ with respect
to the iterations of \textbf{fmincon. }We see that\textbf{\ }$|\nabla
J_{\lambda ,\beta }( v,w,p,q) |$ $\approx 0.01$ after 20 iterations. Then
this value decreases very slowly with iterations. These justify the stopping
criterion (\ref{7.3}). We see on Figure \ref{plot_re_A2_diff_lambda} that
with $\lambda =3$ the reconstructions of both the shape of the inclusion and
the value of the contrast (\ref{8.2040}) are accurate ones.

\begin{figure}[tbph]
\centering
\includegraphics[width = 2.8in]{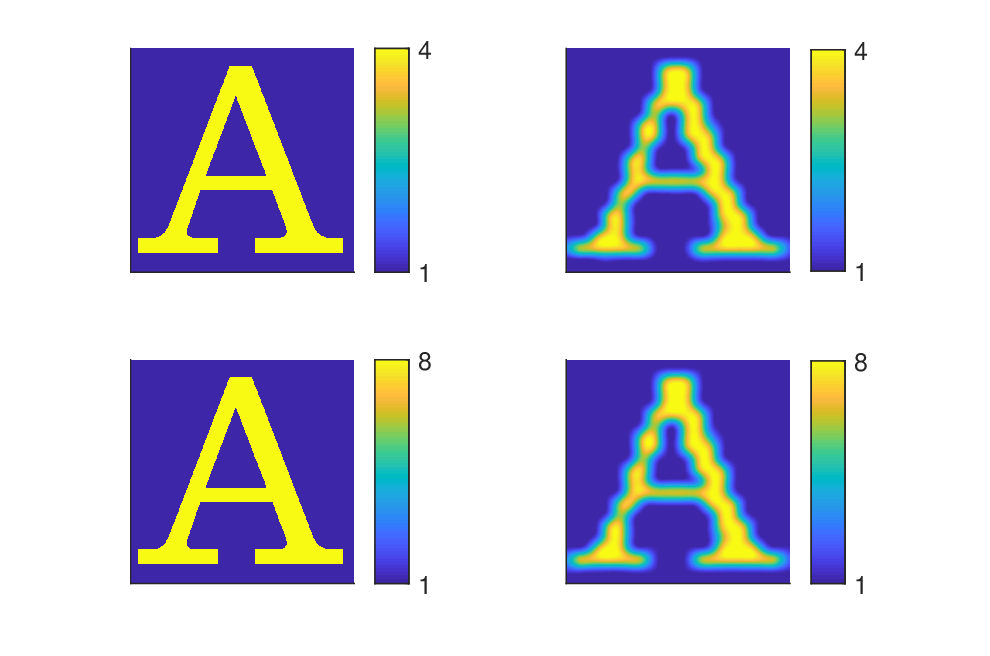}
\caption{Test 2. Exact (left) and reconstructed (right) coefficient $k(x)$
with $c_{a}=4$ (first row) and $c_{a}=8$ (second row) inside of the letter
`A' as in (\protect\ref{8.203}). The inclusion/background contrasts in (%
\protect\ref{8.204}) are respectively $4:1$ and $8:1$. The reconstructions
of both shapes of inclusions and the inclusion/background contrasts (\protect
\ref{8.2040}) are accurate.}
\label{plot_re_A_4_8}
\end{figure}

\textbf{Test 2.} We test the case when the inclusion in (\ref{8.203}) has
the shape of the letter `$A$' for different values of the parameter $%
c_{a}=4,8$ inside of the letter `$A$'. Hence, by (\ref{8.204}) the
inclusion/background contrasts now are respectively $4:1$ and $8:1$.
Computational results are displayed on Figure \ref{plot_re_A_4_8}. One can
observe that shapes of inclusions are imaged accurately. In addition, the
computed inclusion/background contrasts (\ref{8.2040}) are accurate. 

\begin{figure}[tbph]
\centering
\includegraphics[width = 3in]{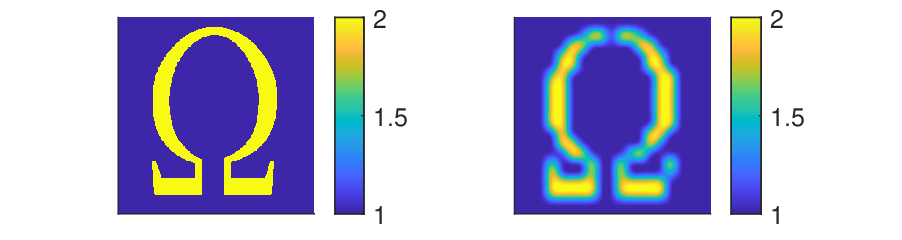}
\caption{Test 3. Exact (left) and reconstructed (right) coefficient $k(x)$,
where the function $k(x)$ is given in \eqref{8.203} with $c_{a}=2$ inside of
the letter `$\Omega $'. The reconstructions of both the shape of the
inclusions and the inclusion/background contrast (\protect\ref{8.2040}) are
accurate.}
\label{plot_re_Omega}
\end{figure}

\textbf{Test 3.} We test the case when the coefficient $k(x)$ in (\ref{8.203}%
) has the shape of the letter `$\Omega $' with $c_{a}=2$ inside of it.
Results are presented on Figure \ref{plot_re_Omega}. We again observe an
accurate reconstruction. 

\begin{figure}[tbph]
\centering
\includegraphics[width = 3in]{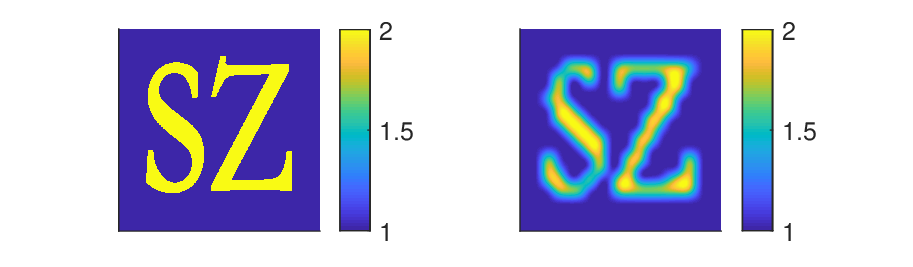}
\caption{Test 4. Exact (left) and reconstructed (right) coefficient $k(x)$,
where the function $k(x)$ is given in \eqref{8.203} with $c_{a}=2$ inside of
two letters `SZ'. The reconstruction is worse than the one for the case of
the single letter `$\Omega $' in Figure \protect\ref{plot_re_Omega}.
Nevertheless, the reconstructions of shapes of both letters are still
accurate. In addition, the computed inclusion/background contrasts in 
\eqref{8.2040} are accurately reconstructed in both letters.}
\label{plot_re_SZ}
\end{figure}

\textbf{Test 4.} We test the reconstruction for the case when the inclusion
in (\ref{8.203}) has the shape of two letters `SZ' with $c_{a}=2$ in each of
them. S and Z are two letters in the name of the city (Shenzhen) were the
second author resides. The results are exhibited on Figure \ref{plot_re_SZ}.
Reconstructions of shapes of both letters as well as of the computed
inclusions/background contrasts (\ref{8.2040}) in both letters are accurate
ones. 

\textbf{Test 5.} We consider the case when the random noisy is present in
the data in (\ref{8.210}) with $\delta =0.03$ and $\delta =0.05$, i.e. with
3\% and 5\% noise level. We test the reconstruction for the cases when the
inclusion in (\ref{8.203}) has the shape of either the letter `$A$' or the
letter `$\Omega $' with $c_{a}=2$. The results are displayed on Figure \ref%
{plot_re_A2_AddNoiseOn_boundary}. One can observe accurate reconstructions
in all four cases. In particular, the inclusion/background contrasts in (\ref%
{8.2040}) are reconstructed accurately.

\begin{figure}[tbph]
\centering
\includegraphics[width = 4.5in]{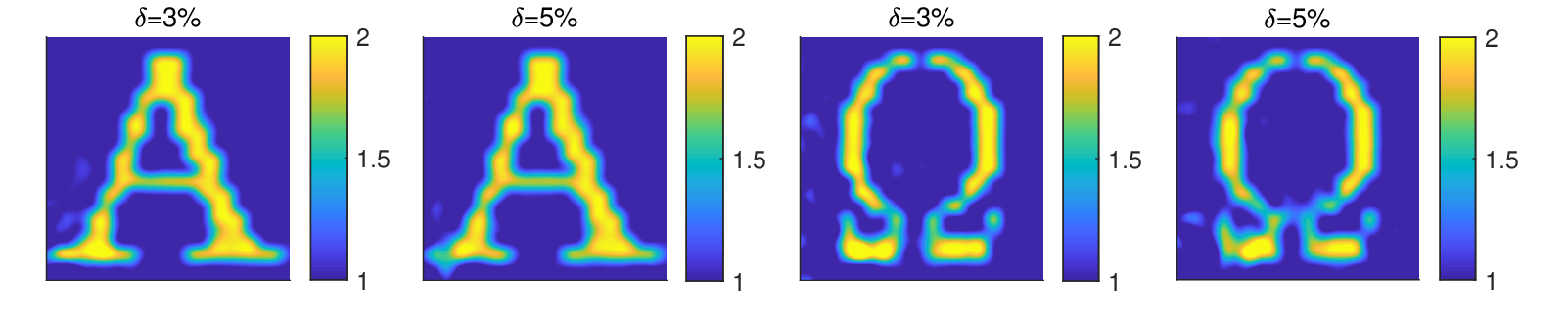}
\caption{Test 5. Reconstructed coefficient $k(x)$ with the shape of letters
`A' or the letter `$\Omega $' with $c_{a}=2$ from noisy data \eqref{8.210}
with $\protect\delta =0.03$ and $\protect\delta =0.05$, i.e. with 3\% and
5\% noise level. In all these four cases, reconstructions of both shapes of
inclusions and inclusion/background contrasts in \eqref{8.2040} are accurate.
}
\label{plot_re_A2_AddNoiseOn_boundary}
\end{figure}


\end{document}